\documentclass{llncs}

\usepackage{amsmath}
\usepackage{amsfonts}
\usepackage{amssymb}
\usepackage{array}
\usepackage{latexsym}
\usepackage{exscale,cmmib57}
\usepackage{llncsdoc}
\usepackage{graphicx}
\usepackage{url}

\def \bx {\mathbf{x}}

\newcommand*\Bell{\ensuremath{\boldsymbol\ell}}

\begin{document}

\title{Multilevel sparse grids collocation for linear partial differential equations, with tensor product smooth basis functions}
%
%
\author{Yangzhang Zhao, Qi Zhang, Jeremy Levesley}
%
%
%
\institute{ Department of Mathematics, University of Leicester, \\ University Road, Leicester, LE1 7RH, United Kingdom\\
{\tt \{ yz177,\ qz49,\ jl1 \} @le.ac.uk} }

\maketitle              

\begin{abstract}

Radial basis functions have become a popular tool for approximation and solution of partial differential equations (PDEs). The recently proposed multilevel sparse interpolation with kernels (MuSIK) algorithm proposed in \cite{Georgoulis} shows good convergence. In this paper we use a sparse kernel basis for the solution of PDEs by collocation. We will use the form of approximation proposed and developed by Kansa \cite{Kansa1986}. We will give numerical examples using a tensor product basis with the multiquadric (MQ) and Gaussian basis functions. This paper is novel in that we consider space-time PDEs in four dimensions using an easy-to-implement algorithm, with smooth approximations. The accuracy observed numerically is as good, with respect to the number of data points used, as other methods in the literature; see \cite{Langer1,Wang1}.

\end{abstract}

\section{Introduction}

During the last few decades since radial basis functions (RBFs) were proposed by Hardy \cite{Hardy} for numerical approximation, they have been applied to a wide range of applications from mathematics, geophysics, physics to engineering and finance. In this paper we will use tensor products of the infinitely differentiable univariate functions
\begin{eqnarray*}
Multiquadric &:& \phi_{c}(x) = \sqrt{x^{2} + c^{2}},\\
Gaussian &:& \psi_{c} (x) = e^{- \frac{x^{2}}{c^{2}}}.
\end{eqnarray*}
The basis function for approximation is then of the form
\begin{equation*}
\Phi_{ \mathbf{c} } \left( \mathbf{x} \right) = \prod^{d}_{i=1} \mu_{c_{i}} \left( x_{i} \right),
\end{equation*}
where $\mu$ is either $\phi$ or $\psi$. This is not strictly speaking RBF approximation in general, though for the Gaussian basis function, since
\begin{equation*}
\prod^{d}_{i=1} \exp \left( -x^{2}_{i} \right) = \exp \left( - \left( \sum^{d}_{i=1} x^{2}_{i} \right) \right) = \exp \left( - \| \mathbf{x} \|^{2} \right),
\end{equation*}
we obtain a univariate function of the norm (an RBF).

\begin{figure}[!htb]
\begin{center}
   \includegraphics[width=6cm]{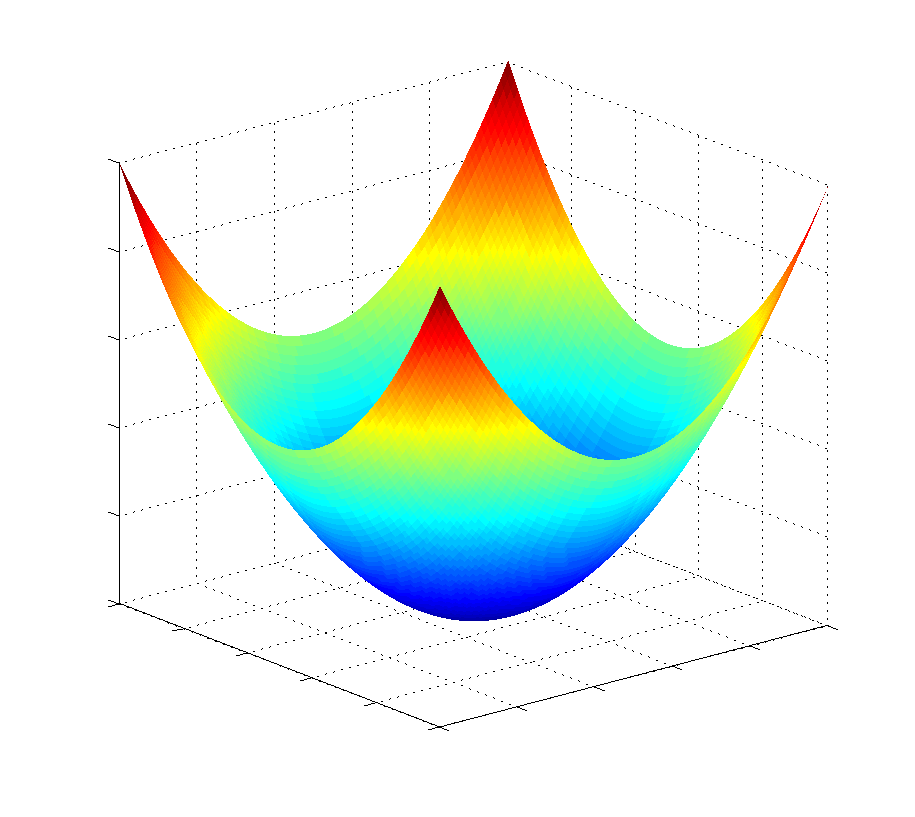}
   \includegraphics[width=6cm]{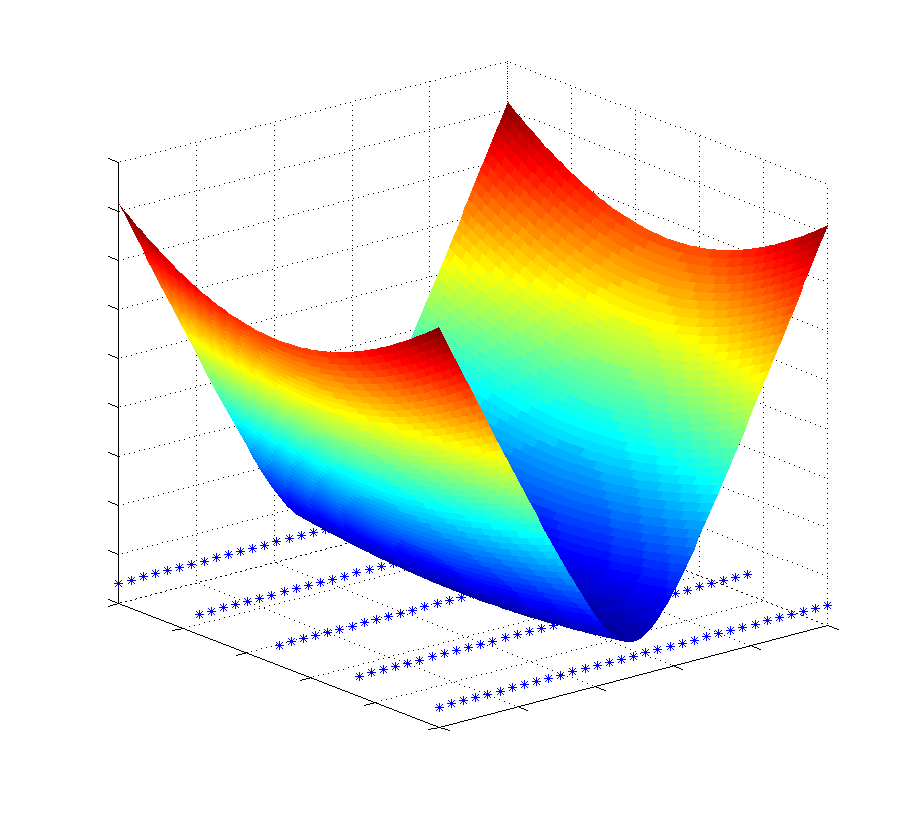}
\end{center}
\caption{An example of normal MQ and anisotropic tensor MQ funcions in 2D. The anisotropic function to the right is scaled appropriately for the inisotropic grid shown.} \label{Fig:ARBF_MQ}
\end{figure}

\begin{figure}[!htb]
\begin{center}
   \includegraphics[width=6cm]{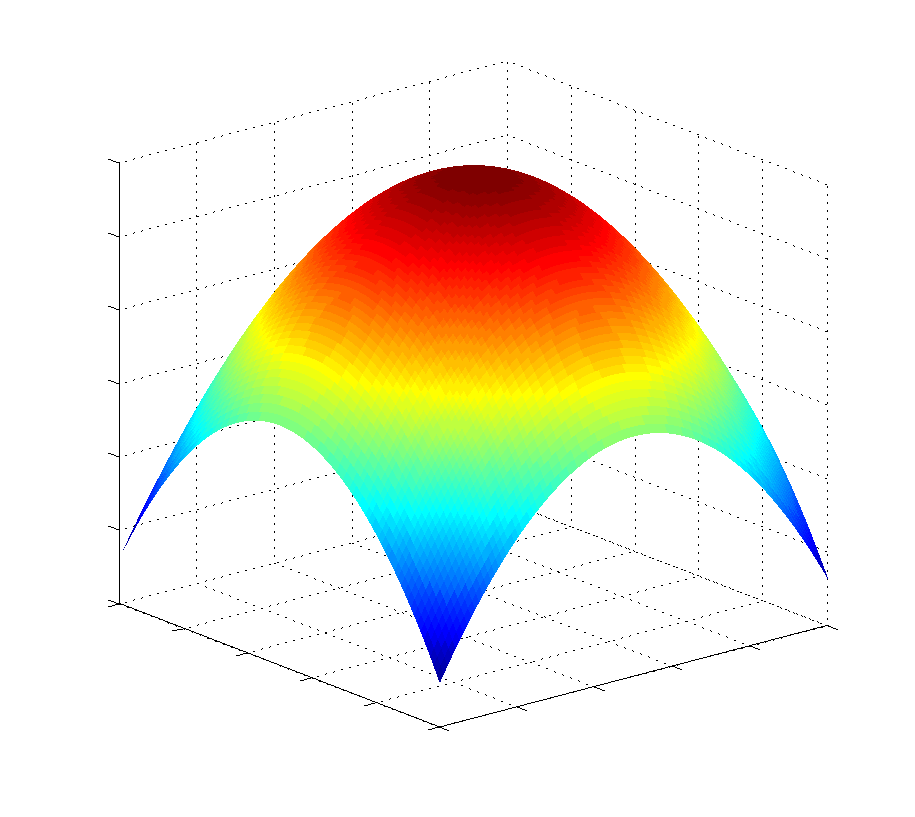}
   \includegraphics[width=6cm]{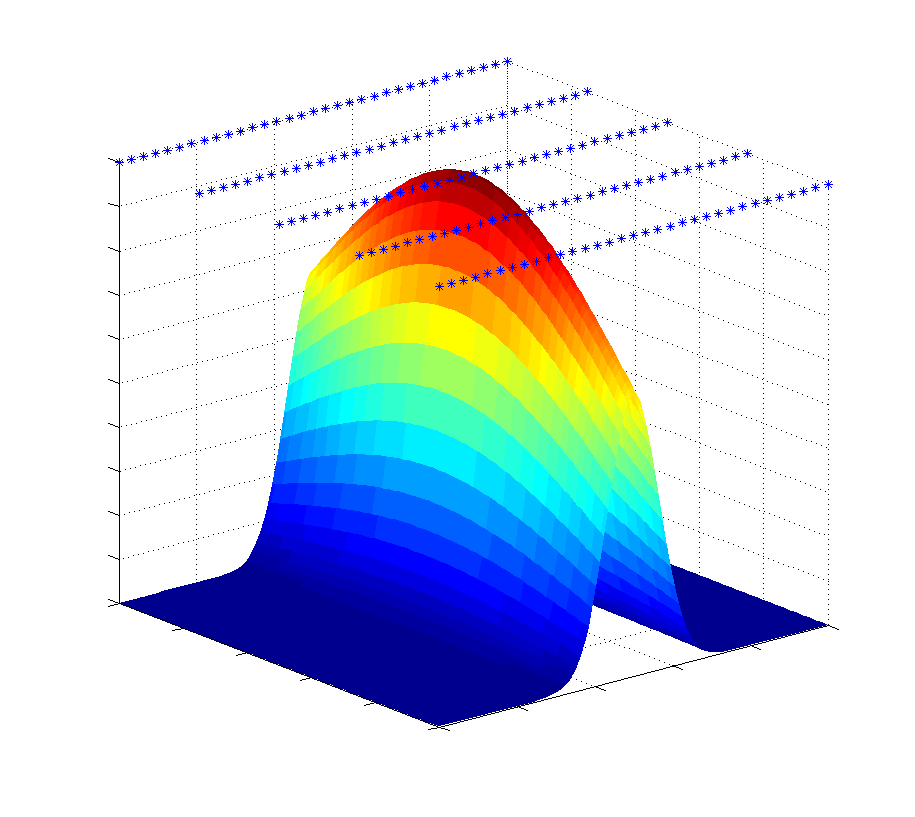}
\end{center}
\caption{An example of normal Gaussian and anisotropic tensor Gaussian functions in 2D. The anisotropic function to the right is scaled appropriately for the inisotropic grid shown.} \label{Fig:ARBF_Gauss}
\end{figure}

In the definition of the multiquadric and Gaussian there is a parameter $c$
which we call the shape parameter. This is used to scale the approximation
basis in various directions depending on the resolution of the data points in that direction; see Figures \ref{Fig:ARBF_MQ} and \ref{Fig:ARBF_Gauss}. In the first we plot $\Phi_{\left[1,1\right]}\left(\mathbf{x} \right)=\phi_{1}(x_{1})\phi_{1}(x_{2})$ and $\Phi_{\left[1/2,1/32\right]}\left(\mathbf{x} \right)=\phi_{1/2}(x_{1})\phi_{1/32}(x_{2})$, and in the second $\Phi_{\left[1,1\right]}\left(\mathbf{x} \right)=\psi_{1}(x_{1})\psi_{1}(x_{2})$ and $\Phi_{\left[1/2,1/32\right]}\left(\mathbf{x} \right)=\psi_{1/2}(x_{1})\psi_{1/32}(x_{2})$. The scaling matches the anisotropic grid shown next to the surface. We call the basis functions with different shape in each direction anisotropic basis functions.

More recently RBFs have been employed in the solution of PDEs \cite{Babuska,Brown,Farrell,Fasshauer,Fornberg,Franke1,Franke2,Giesl,Hales,Nobile,Schaback}. Suppose our PDE is
\begin{eqnarray}
\mathcal{L} u &=& f \quad \text{in} \ \Omega, \label{PDE1} \\
u &=& g \quad \text{on} \ \partial \Omega. \label{PDE2}
\end{eqnarray}

There are two distinct collocations methods using RBFs in this context, termed symmetric and non-symmetric collocation. The latter was introduced by Kansa \cite{Kansa1,Kansa2,Kansa1986} and involves the expansion of the solution of the PDE in a combination of RBFs:
\begin{equation*}
\tilde{u} \left( \mathbf{y} \right) = \sum_{i} \alpha_{i} \Phi_{\mathbf{c}_{i}} \left( \mathbf{y} - \mathbf{x}_{i} \right),
\end{equation*}
where the nodes $\mathbf{x}_{i} \in \Omega \cup \partial \Omega$. The PDE is applied to this expansion and collocation is used to compute coefficients in the expansion. Currently there is no proof that this method is stable in the sense that the collocation system is invertible. However, the method remains simple to implement and shows good convergence. In symmetric collocation, developed by \cite{Fasshauer1}, the solution of the PDE is written in the form:
\begin{equation*}
\tilde{u} \left( \mathbf{y} \right) = \sum_{i} \alpha_{i} \mathcal{L} \Phi_{\mathbf{c}_{i}} \left( \mathbf{y} - \mathbf{x}_{i} \right) + 
\sum_{j} \beta_{j} \Phi_{\mathbf{c}_{j}} \left( \mathbf{y} - \mathbf{z}_{j} \right),
\end{equation*}
where now $\mathbf{x}_{i} \in \Omega$, and $\mathbf{z}_{j} \in \partial \Omega$. Now the collocation system which arises is symmetric and for specific choices of positive definite RBFs (the Gaussian for instance) the system can be proven to be invertible (see \cite{Fasshauer1}).

Due to the simplicity of implementation we will use non-symmetric collocation in this paper. We will explore the use of symmetric collocation with sparse grids in a follow-up article.

One of the advantages in using radial basis functions is the ease of implementation
in high dimensional problems, though this is of no practical consequence if we cannot mitigate the so-called \emph{curse of dimensionality}. The sparse grid methodology which will be described later is our chosen route to doing this, though this of itself is not new; see e.g. \cite{Bungartz,Ganapathysubramanian,Griebel,Nobile}. If the discretisation parameter (the smallest separation distance between two points) is $h$, then in a full grid method the number of points is of the order $h^{-d}$. In the sparse grid algorithm we have a number of points of the order $h^{-1} | \log h |^{d-1}$.

A new feature in this paper is that we use smooth kernels in the sparse grid
algorithm. This means that there is no restriction to the potential convergence rates we may get related to approximations which are of finite smoothness. For instance, for univariate linear or cubic B-spline approximation we should expect no better than $\mathcal{O} \left( h^{2} \right)$ or $\mathcal{O} \left( h^{4} \right)$ rate of convergence respectively, regardless of the smoothness of the function we are approximating. On the other hand, trigonometric approximation of periodic analytic functions converges with exponential rates since the approximating functions are also analytic; see e.g. \cite{rivlin}. As we will see in our examples in Section 5, the rate of convergence of the methods depend on the choice of the shape parameter, and that on the loglog plots the rate curves are convex, suggesting that we are achieving rates faster than any polynomial.

We compare with Wang et al. \cite{Wang1}, who use a sparse grid algorithm with piecewise polynomials. It is usual to treat space and time separately, but we treat them together, as do Schwab and Stevenson in \cite{schwab}, and we compare our algorithm with Langer et al. \cite{Langer1}. The reason we treat them together is that in a method separating space and time we would expect to multiply the complexity of the space part of the algorithm by the number of time steps. Use of the sparse grid algorithm leads to multiplication by the log of the time discretisation.

An advantage of our algorithm is that it is relatively straightforward to code. A disadvantage as will be seen in the numerical results is that we suffer ill-conditioning problems as the shape parameters increase. This is an instance of the well-known uncertainty principle in RBF approximation \cite{SchabackU}. As the approximating basis function gets smoother (the shape parameter increases), the approximation of smooth functions gets better, but the conditioning of the approximation equations gets worse. In future work we will seek to mitigate the ill-conditioning while maintaining the approximation power.

Of course, there are other methods that have been applied to solving PDEs
using RBFs. In particular, we should point out other multiscale or multilevel methods such as detailed in \cite{Farrell,Fasshauer}. Related to more traditional methods we have RBF finite difference methods \cite{Fornberg1}, RBF finite element methods \cite{Hu}, and RBF partition of unity methods \cite{Shcherbakov}.

The collocation method described in the sequel is based on the recently developed approximation method called Multilevel Sparse Interpolation with Kernels (MuSIK) \cite{Georgoulis}. This method has been used to solve approximation problems in up to 5 dimensions, and quadrature problems in up to 10 dimensions.

In Section 3 we introduce RBF collocation, and in Section 2 we describe
the implementation of the MuSIK-C, multilevel sparse grid kernel collocation. In Section 5 we apply this algorithm to a number of PDEs, including elliptic and parabolic time-dependent PDEs (the heat equation). In the latter we treat the time as one spatial dimension, so do not do any time-stepping such as in the method of lines. In this paper we aim to demonstrate that our method has potential to solve PDEs in high dimensions. Therefore, our low dimensional examples demonstrate this potential. We do not pretend that we are capable of solving non-smooth PDEs on complicated domains, which are very much the domain of more well-established methods. We restrict ourselves to four dimensional problems (three in space and one in time), but the results that we achieve indicate that higher dimensional problems, to be considered in future work, are tractable.

\section{Sparse grids} \label{sparse}

Multilevel sparse kernel-based interpolation (MuSIK) is described in \cite{Georgoulis}. The collocation method is almost identical. We begin by describing sparse grids.

Let $\Omega = \left[0,1 \right]^{d}$ and $\partial \Omega$ be its boundary. Furthermore, let $\Bell=\left( l_{1},\dots,l_{d} \right) \in \mathbb{N}^{d}$ be a multi-index, and for $0 \leq k_{i} \leq 2^{l_{i}}$, $i=1,2,\dots,d$, $x_{i,k_{i}}=k_{i}2^{-l_{i}}$, be a uniform partition of $\left[0,1 \right]$. Then we define the family of grids which are uniform in each direction $\mathbb{X}_{\Bell} = \{x_{i,k_{i}}, 0 \leq k_{i} \leq 2^{l_{i}}, i=1,2,\dots,d \}$. The total number of nodes in $\mathbb{X}_{\Bell}$ is less than $2^{\| \Bell \|_{1}+1}$, where
\begin{equation*}
\| \Bell \|_{1} = \sum^{d}_{i=1} l_{i}
\end{equation*}
is the one-norm of $\Bell$.

\begin{figure}
\begin{center}
\begin{tabular}{m{1.8cm}m{0.2cm}m{1.6cm}m{0.2cm}m{1.6cm}m{0.2cm}m{1.6cm}m{0.2cm}m{1.6cm}}
   \includegraphics[width=1.9cm]{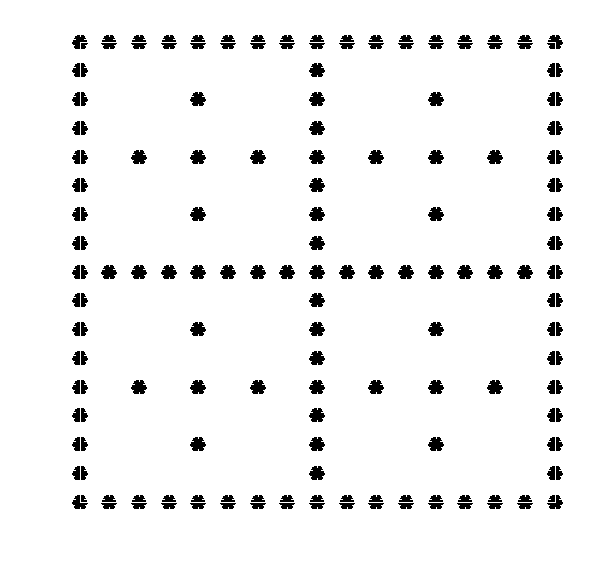} &=
      &\includegraphics[width=1.7cm]{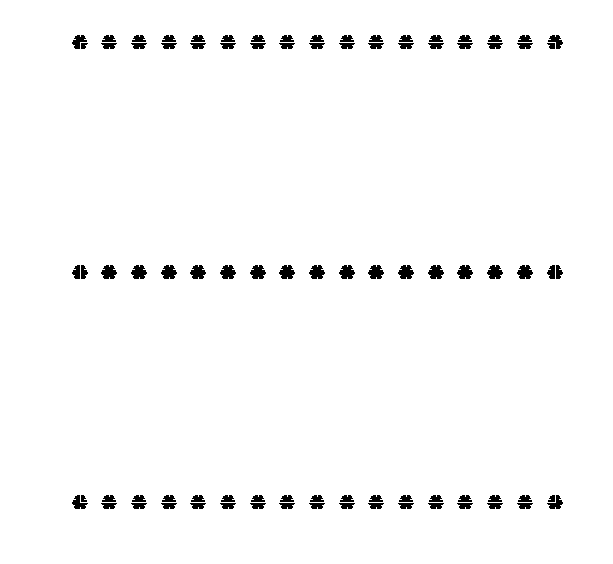} 
     & $\ \cup$  &\includegraphics[width=1.7cm]{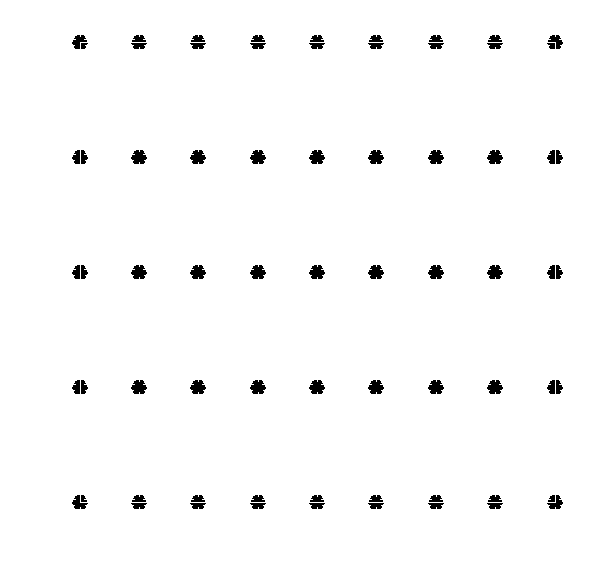} 
     & $\ \cup$ & \includegraphics[width=1.7cm]{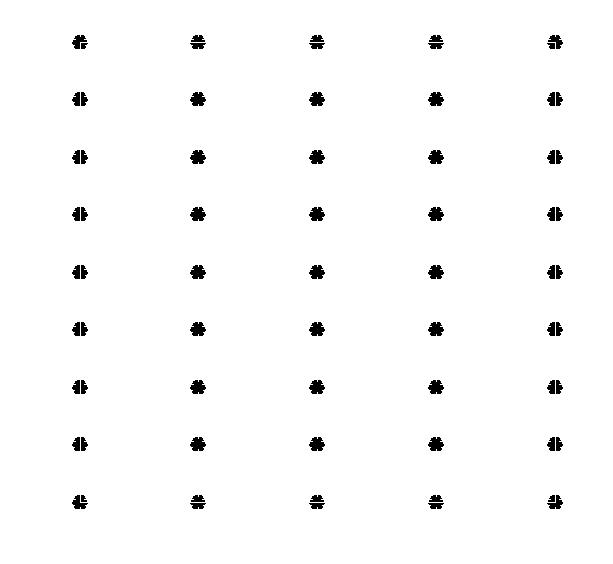} 
     & $\ \cup$ & \includegraphics[width=1.7cm]{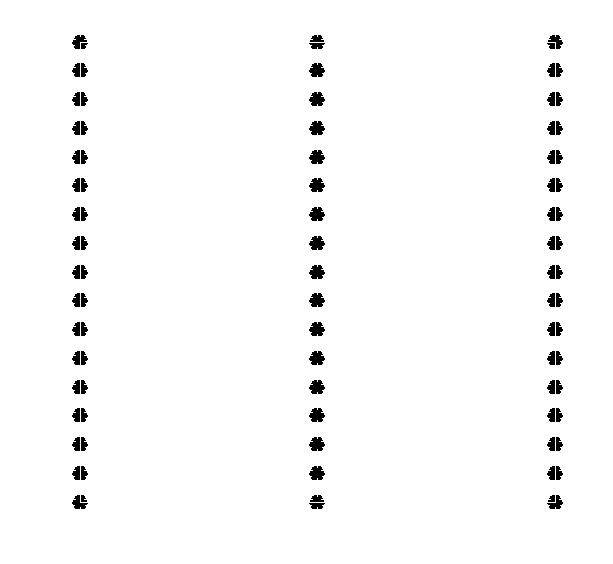}  
\end{tabular}
\end{center}
\caption{Sparse grid $\tilde{\mathbb{X}}^{4,2}$ via \eqref{sparse_grid_def}. } \label{Fig:Sparse grid def}
\end{figure}

\begin{figure}
\begin{center}
\begin{tabular}{m{1.8cm}m{0.2cm}m{1.8cm}m{0.2cm}m{1.8cm}m{0.2cm}m{1.8cm}m{0.2cm}m{1.8cm}}
   \includegraphics[width=1.9cm]{gridS4.png} &=
      &\includegraphics[width=1.8cm]{gridX41.png} 
     & $\oplus$  &\includegraphics[width=1.8cm]{gridX32.png} 
     & $\oplus$ & \includegraphics[width=1.8cm]{gridX23.png} 
     & $\oplus$ & \includegraphics[width=1.8cm]{gridX14.png}  \\
     & $\ominus$  &\includegraphics[width=1.8cm]{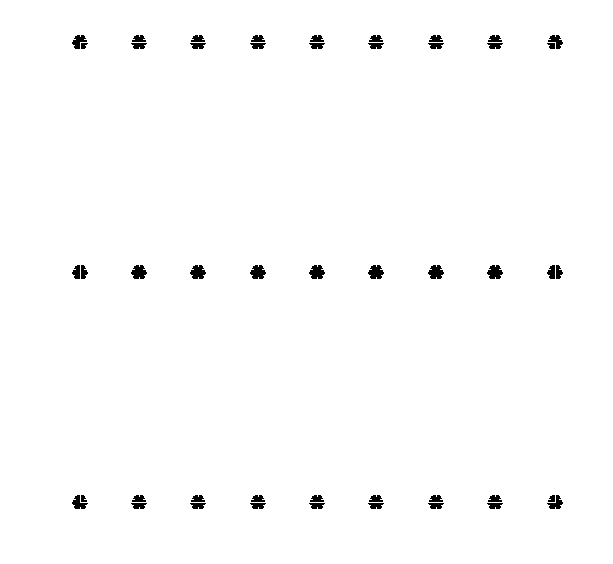} 
     & $\ominus$ & \includegraphics[width=1.8cm]{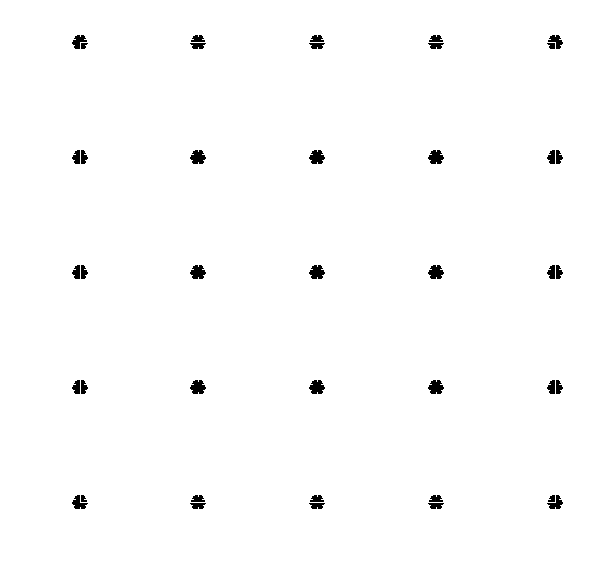} 
     & $\ominus$ & \includegraphics[width=1.8cm]{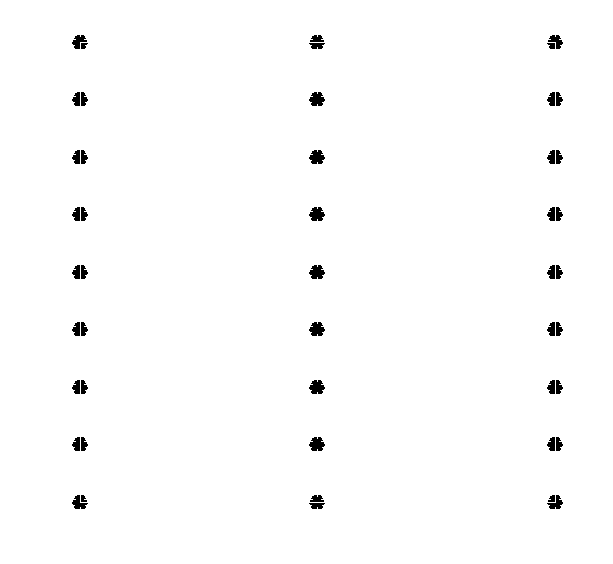}  &
\end{tabular}
\end{center}
\caption{The Boolean sum of the sparse grid $\tilde{\mathbb{X}}^{4,2}$ from full grids at two different levels, via \eqref{sparse combination grid}.} \label{Fig:Sparse grid def2}
\end{figure}

The sparse grid
\begin{equation} \label{sparse_grid_def}
\tilde{\mathbb{X}}^{n,d} := \bigcup_{|\Bell|_{1} = n+(d-1)} \mathbb{X}_{\Bell}.
\end{equation}
See Figure \ref{Fig:Sparse grid def} for an illustration with $n=4$ and $d=2$. A great insight in sparse
grid technology is that this sparse grid can be seen as a Boolean sum of grids at different levels \cite{Delvos}. This means that we can write
\begin{equation} \label{sparse combination grid}
\tilde{\mathbb{X}}^{n,d} := \sum^{d-1}_{q=0} \left( -1 \right)^{q} \binom{d-1}{q} \sum_{\| \Bell \|_{1} = n+(d-1)-q} \mathbb{X}_{\Bell};
\end{equation}
see Figure \ref{Fig:Sparse grid def2}. Here we interpret the plus and minus signs as inclusion or exclusion of points.

The hope is that we can achieve the same approximation power using the sparse grid as we do using the full grid. The number of points in the sparse grid is bounded by
\begin{equation*}
2^{n+d} \left( \# \{ \Bell: \| \Bell \|_{1} = n + (d-1) \} \right) \leq 2^{n+d} \frac{(n+d-1)^{d-1}}{(d-1)!} \leq C(d) 2^{n} n^{d-1},
\end{equation*}
where $C(d)$ is a positive constant which depends on $d$, but not $n$. Since $N=2^{n}+1$ is the number of points along one side of the grid, we see that we have $\mathcal{O} \left( N \left( \log N \right)^{d-1} \right)$ points compared to $\mathcal{O} \left( N^{d} \right)$ for a full grid.

\section{Multilevel sparse grid collocation using tensor product basis functions}

We compute our sparse grid approximation to the solution of the differential
equation by computing a solution separately on each of the grids in the Boolean sum for the sparse grid in \eqref{sparse combination grid}. Let $\mathbb{B}_{\Bell} = \mathbb{X}_{\Bell} \cap \partial \Omega $ be the boundary nodes of $\mathbb{X}_{\Bell}$, and $\mathbb{I}_{\Bell} = \mathbb{X}_{\Bell} - \mathbb{B}_{\Bell}$ be the interior nodes.

With an abuse of the notation introduced above, let us denote by $\Phi_{\Bell}$ the basis function with shape parameters $c_{i}=C/2^{l_{i}}, i=1,2,\dots,d$, for some constant $C$, which we will specify in the examples in Section \ref{Numeric}. The collocation approximation to the solution of the PDE \eqref{PDE1} and \eqref{PDE2} on $\mathbb{X}_{\Bell}$ is
\begin{equation}
u_{\Bell} (\mathbf{y}) = \sum_{\mathbf{x} \in \mathbb{X}_{\Bell}} \Phi_{\Bell} \left( \mathbf{y} - \mathbf{x} \right),
\end{equation}
which satisfies
\begin{eqnarray}
\mathcal{L} u_{\Bell} (\mathbf{x}) &=& f(\mathbf{x}), \quad \mathbf{x} \in \mathbb{I}_{\Bell}, \\
u_{\Bell}(\mathbf{x}) &=& g(\mathbf{x}), \quad \mathbf{x} \in \mathbb{B}_{\Bell}.
\end{eqnarray}

Mimicking the Boolean decomposition of the sparse grid into uniform grids, we form the sparse grid approximation to $u$
\begin{equation}
u^{n,d} := \sum^{d-1}_{q=0} \left( -1 \right)^{q} \binom{d-1}{q} \sum_{\| \Bell \|_{1} = n+(d-1)-q} u_{\Bell},
\end{equation}
by combining the solutions on the uniform subgrids. This is called the combination technique, introduced by Delvos \cite{Delvos}. We call this algorithm SIK-C.

We will see in the numerical examples below, and as is also observed in
interpolation \cite{Georgoulis} and quasi-interpolation \cite{Usta}, that $u^{n,d} \nrightarrow u$ as $n \rightarrow \infty$ for the Gaussian basis function, and convergence is slow for the multiquadric. This is because we are scaling the shape parameter at each level of approximation, and convergence results are really only available for fixed shape parameter. In order to obtain convergence we employ a multilevel refinement strategy, exactly as in \cite{Georgoulis}.

The \emph{multilevel sparse kernel-based collocation} (MuSIK-C, for short) algorithm is initialised by computing the SIK-C solution $u^{n_{0},d}$ on the coarsest sparse grid $\tilde{\mathbb{X}}^{n_{0},d}$ and set $\Delta_{0}:=u^{n_{0},d}$. Then, for $k=1,2,\dots$, $\Delta_{k}$ is the SIK-C solution to the residual $\Gamma_{k}:=u - \sum_{j=0}^{k-1} \Delta_{j}$ on $\tilde{\mathbb{X}}^{n_{0}+k,d}$. The resulting MuSIK-C solution is then given by
\begin{equation*}
u^{n,d}_{\rm ML} := \sum^{n}_{j=0} \Delta_{j}.
\end{equation*}

\section{Extrapolation}

Extrapolation is a process by which we can accelerate convergence if we understand how the asymptotic error behaves. The most well-known instance of the process is Romberg integration, a description which can be found in e.g. \cite[Page 211]{Burden}. What we will show is that, even if we do not know the precise form of the asymptotic expansion, we can reduce the error significantly using a simple extrapolation process. The coefficient in the extrapolation is estimated using the most recent errors. In this implementation we do not use the accelerated solutions in any further acceleration process.

We estimate the convergence rate for our acceleration from the errors of the previous two steps:
\begin{equation*}
\beta \approx - \frac{\log E_{N_{2}}(f) - \log E_{N_{1}}(f) }{\log N_{2} - \log N_{1} }.
\end{equation*}
In the following numerical examples we will see that this estimate for $\beta$ provides significant improvement in the error, though we do not improve the rate of approximation. For the latter we would need to know the rate exactly.

\section{Numerical experiments} \label{Numeric}

In this section, we will employ MuSIK-C to solve a variety of elliptic and
parabolic partial differential equations in up to four dimensions. In the following tables, "Cond" represents the condition number at level $n$, "Nodes" is the number of sparse grid centers used, and we measure the errors in SIK-C and MuSIK-C respectively with
\begin{eqnarray*}
E^{n,d}_{\rm SIK-C} (\bx) &=& \max_{\bx \in \mathbf{T}} \left| u(\bx)- u^{n,d} (\bx) \right|, \\
E^{n,d}_{\rm MuSIK-C} (\bx) &=& \max_{\bx \in \mathbf{T}} \left| u(\bx)- u_{\rm ML}^{n,d} (\bx) \right|,
\end{eqnarray*}
where $\mathbf{T}$ refers to a test set of points which will be made explicit in each example. We will also use both multiquadric and Gaussian basis functions for comparison. We will state the value of the parameter $C$ referred to at the start of Section 3. Correspondingly, we define the slope $\rho$ for two adjacent points in different methods, for instance
\begin{eqnarray*}
\rho^{n+1,d}_{\rm MuSIK-C} = \frac{ \log \left( E_{\rm MuSIK-C}^{n+1,d} \right) - \log \left( E_{\rm MuSIK-C}^{n,d} \right) }{ \log \left( {\rm Nodes}^{n+1,d} \right) - \log \left( {\rm Nodes}^{n,d} \right)}.
\end{eqnarray*}
Here ${\rm Nodes}^{n,d}$ means the number of nodes at the $nth$ level in $d$ dimensions. To improve the performance of MuSIK-C further, we employ extrapolation and measure the error and the slope as
\begin{eqnarray*}
E_{\rm extra}^{n+1,d} &=& \max_{x \in \mathbf{T}} \left| u(\mathbf{x}) - u^{n+1,d}_{\rm extra} \right|, \\
\rho^{n+2,d}_{\rm extra} &=& \frac{ \log \left( E_{\rm extra}^{n+2,d} \right) - \log \left( E_{\rm extra}^{n+1,d} \right) }{ \log \left( {\rm Nodes}^{n+2,d} \right) - \log \left( {\rm Nodes}^{n+1,d} \right)}.
\end{eqnarray*}
In the above $\beta^{n+1,d} = - \rho^{n+1,d}_{\rm MuSIK-C} $ and
\begin{eqnarray*}
u^{n+1,d}_{\rm extra} &=& \frac{ \left( \frac{\rm Nodes^{n+1,d}}{\rm Nodes^{n,d}} \right)^{\beta^{n+1,d}} u^{n+1,d}_{\rm MuSIK-C} - u^{n,d}_{\rm MuSIK-C}}{\left( \frac{\rm Nodes^{n+1,d}}{\rm Nodes^{n,d}} \right)^{\beta^{n+1,d}}-1}.
\end{eqnarray*}

In terms of convergence rate with regard to point spacing, since subsequent levels have half the point spacing we can compute the rate as
$$
\rho_h=\log_2 \left ( {E_{\rm MuSIK-C}^{n,d} \over E_{\rm MuSIK-C}^{n+1,d}} \right )
$$
We report these numbers in each example.

\subsection{Elliptic examples}

\begin{example} \label{EE 2D 1}
In this example, we solve the following two-dimensional problem on $\Omega = \left( 0,1 \right)^{2}$
\begin{equation}
\Delta u(\mathbf{x}) = - \pi^{2} \sin(\pi x_{1} x_{2}) (x_{1}^{2} + x_{2}^{2}), \quad \mathbf{x} \in \Omega,
\end{equation}
with boundary conditions
\begin{equation}
u(\mathbf{x}) = \sin(\pi x_{1} x_{2}), \quad \mathbf{x} \in \partial \Omega.
\end{equation}
The exact solution is a two-dimensional non-tensor product function
\begin{equation}
u(\mathbf{x}) = \sin(\pi x_{1} x_{2}).
\end{equation}
\end{example}

\begin{table}[!htb]
\begin{center}
\addtolength{\tabcolsep}{-1pt}
\begin{tabular}{||c|c|c|c|c|c|c||}
  \hline
\ Level \ & \ Nodes \ & \ Cond \ & \ $E_{\rm MuSIK-C}$ \ & \ $\rho_{\rm MuSIK-C}$ \ & \ $E_{\rm SIK-C}$ \  & \ $\rho_{\rm SIK-C}$ \  \\
\hline \hline
 2  &  21     &1e5  & 4.51e-2 & ---    & 4.51e-2 & --- \\
 3  &  49     &7e5  & 1.61e-2 & -1.21  & 1.69e-2 & -1.16 \\
 4  &  113    &8e6  & 3.85e-3 & -1.71  & 9.55e-3 & -0.68 \\
 5  &  257    &5e7  & 8.66e-4 & -1.82  & 5.91e-3 & -0.58 \\
 6  &  577    &2e8  & 2.09e-4 & -1.76  & 4.15e-3 & -0.44 \\
 7  &  1281   &9e8  & 5.17e-5 & -1.75  & 2.93e-3 & -0.44 \\
 8  &  2817   &4e9  & 1.27e-5 & -1.78  & 2.11e-3 & -0.42 \\
 9  &  6145   &2e10 & 3.13e-6 & -1.80  & 1.50e-3 & -0.43 \\
 10 &  13313  &6e10 & 7.68e-7 & -1.81  & 1.08e-3 & -0.43 \\
 11 &  28673  &3e11 & 1.88e-7 & -1.83  & 7.69e-4 & -0.44 \\
 12 &  61441  &1e12 & 4.62e-8 & -1.84  & 5.27e-4 & -0.50 \\[0.5ex]
 \hline
\end{tabular}
\end{center}
\caption{The performance of multilevel sparse collocation and sparse collocation using the MQ for Example \ref{EE 2D 1} with $C=2$. Max error evaluated at 64,000 Halton points in the whole domain. $\rho_h \approx 2.0$.} \label{Table EE 2D 1 C=2 MQ MuSIK}
\end{table}

\begin{table}[!htb]
\begin{center}
\addtolength{\tabcolsep}{-1pt}
\begin{tabular}{||c|c|c|c|c|c|c||}
  \hline
\ Level \ & \ Nodes \ & \ Cond \ & \ $E_{\rm MuSIK-C}$ \ & \ $\rho_{\rm MuSIK-C}$ \ & \ $E_{\rm SIK-C}$ \  & \ $\rho_{\rm SIK-C}$ \  \\
\hline \hline
 2  &  21     &3e6  & 2.78e-2 & ---    & 2.78e-2 & --- \\
 3  &  49     &6e7  & 7.74e-3 & -1.51  & 6.58e-3 & -1.70 \\
 4  &  113    &6e8  & 1.15e-3 & -2.28  & 2.66e-3 & -1.08 \\
 5  &  257    &2e9  & 2.03e-4 & -2.11  & 1.85e-3 & -0.45 \\
 6  &  577    &2e10 & 4.02e-5 & -2.00  & 1.43e-3 & -0.32 \\
 7  &  1281   &1e11 & 8.59e-6 & -1.93  & 1.07e-3 & -0.36 \\
 8  &  2817   &5e11 & 1.83e-6 & -1.96  & 7.82e-4 & -0.40 \\
 9  &  6145   &2e12 & 3.73e-7 & -2.04  & 5.76e-4 & -0.39 \\
 10 &  13313  &9e12 & 7.52e-8 & -2.07  & 4.16e-4 & -0.42 \\
 11 &  28673  &4e13 & 1.54e-8 & -2.07  & 2.98e-4 & -0.43 \\
 12 &  61441  &1e14 & 9.80e-9 & -0.59  & 2.13e-4 & -0.44 \\[0.5ex]
 \hline
\end{tabular}
\end{center}
\caption{The performance of multilevel sparse collocation and sparse collocation using MQ for Example \ref{EE 2D 1} when $C=3$. Max error evaluated at 64,000 Halton points in the whole domain. $\rho_h \approx 2.3$.} \label{Table EE 2D 1 C=3 MQ MuSIK}
\end{table}

\begin{table}[!htb]
\begin{center}
\addtolength{\tabcolsep}{-1pt}
\begin{tabular}{||c|c|c|c|c||}
  \hline
\ Level \ & \ $E_{\rm extra} (C=2)$  \ & \ $\rho_{\rm extra} (C=2)$ \ & \ $E_{\rm extra} (C=3)$ \  & \ $\rho_{\rm extra} (C=3)$ \  \\
\hline \hline
 3    & 7.41e-3 & ---    & 5.38e-3 & ---   \\
 4    & 2.09e-3 & -1.52  & 6.92e-4 & -2.45 \\
 5    & 3.58e-4 & -2.15  & 8.56e-5 & -2.54 \\
 6    & 8.18e-5 & -1.82  & 1.59e-5 & -2.08 \\
 7    & 2.00e-5 & -1.77  & 3.63e-6 & -1.85 \\
 8    & 4.95e-6 & -1.77  & 6.61e-7 & -2.16 \\
 9    & 1.22e-6 & -1.80  & 1.09e-7 & -2.31 \\
 10   & 2.99e-7 & -1.82  & 2.10e-8 & -2.14 \\
 11   & 7.39e-8 & -1.82  & 4.71e-9 & -1.95 \\
 12   & 1.58e-8 & -2.02  & 2.75e-8 &  2.32 \\[0.5ex]
 \hline
\end{tabular}
\end{center}
\caption{Extrapolations from the multilevel sparse collocation using MQ for Example \ref{EE 2D 1} with differnet constants: $C=2$ and $C=3$. Max error evaluated at 64,000 Halton points in the whole domain. $\rho_h \approx 2.0$.} \label{Table EE 2D 1 MQ MuSIK extra}
\end{table}

From Tables \ref{Table EE 2D 1 C=2 MQ MuSIK} and \ref{Table EE 2D 1 C=3 MQ MuSIK}, it is obvious that the multilevel method really offers advantages in the solutions. The difference between these two tables is the value of the parameter $C$, ($C=2$ and $3$). When $C$ is bigger, we have a smoother basis function. Thus, as we suggested in the introduction, the condition number grows faster, but the convergence rate is faster and solutions are more accurate. At Level 12, in Table 2 we see that a condition number of $1e14$ is effecting the answer adversely. From Table \ref{Table EE 2D 1 MQ MuSIK extra}, we observe that extrapolation reduces the error by a factor of 3 or 4. 

\begin{table}[!htb]
\begin{center}
\addtolength{\tabcolsep}{-1pt}
\begin{tabular}{||c|c|c|c|c|c|c||}
  \hline
\ Level \ & \ Nodes \ & \ Cond \ & \  $E_{\rm MuSIK-C}$ \ & \ $\rho_{\rm MuSIK-C}$ \ & \ $E_{\rm SIK-C}$ \  & \ $\rho_{\rm SIK-C}$ \  \\
\hline \hline
 2  &  21     &2e4  & 2.82e-2 & ---    & 2.82e-2 & --- \\
 3  &  49     &8e4  & 1.69e-2 & -0.60  & 2.17e-2 & -0.31 \\
 4  &  113    &2e6  & 4.44e-3 & -1.60  & 2.32e-2 &  0.08 \\
 5  &  257    &3e7  & 9.69e-4 & -1.85  & 2.43e-2 &  0.06 \\
 6  &  577    &2e8  & 2.43e-4 & -1.71  & 2.48e-2 &  0.02 \\
 7  &  1281   &8e8  & 6.21e-5 & -1.72  & 2.51e-2 &  0.02 \\
 8  &  2817   &3e9  & 1.57e-5 & -1.74  & 2.52e-2 &  0.01 \\
 9  &  6145   &1e10 & 3.93e-6 & -1.78  & 2.52e-2 & -0.00 \\
 10 &  13313  &5e10 & 9.88e-7 & -1.79  & 2.47e-2 & -0.03 \\
 11 &  28673  &2e11 & 2.44e-7 & -1.82  & 2.52e-2 &  0.03 \\
 12 &  61441  &9e11 & 5.92e-8 & -1.86  & 2.43e-2 & -0.05 \\[0.5ex]
 \hline
\end{tabular}
\end{center}
\caption{Multilevel sparse collocation compared with sparse collocation using the Gaussian for Example \ref{EE 2D 1} when $C=2$. Max error evaluated at 64,000 Halton points in the whole domain. $\rho_h \approx 2.9$.} \label{Table EE 2D 1 C=2 G MuSIK}
\end{table}

\begin{table}[!htb]
\begin{center}
\addtolength{\tabcolsep}{-1pt}
\begin{tabular}{||c|c|c|c|c|c|c||}
  \hline
\ Level \ & \ Nodes \ & \ Cond \ & \  $E_{\rm MuSIK-C}$ \ & \ $\rho_{\rm MuSIK-C}$ \ & \ $E_{\rm SIK-C}$ \  & \ $\rho_{\rm SIK-C}$ \  \\
\hline \hline
 2  &  21     &9e5  & 2.31e-2 & ---    & 2.31e-2 & --- \\
 3  &  49     &3e8  & 5.98e-3 & -1.60  & 5.87e-3 & -1.62 \\
 4  &  113    &4e9  & 4.75e-4 & -3.03  & 2.74e-3 & -0.91 \\
 5  &  257    &5e10 & 5.51e-5 & -2.62  & 2.36e-3 & -0.18 \\
 6  &  577    &2e13 & 8.02e-6 & -2.38  & 2.35e-3 & -0.00 \\
 7  &  1281   &2e15 & 1.12e-6 & -2.47  & 2.36e-3 &  0.01 \\
 8  &  2817   &2e16 & 1.50e-7 & -2.55  & 2.37e-3 &  0.01 \\
 9  &  6145   &3e17 & 1.99e-8 & -2.59  & 2.37e-3 & -0.00 \\
 10 &  13313  &3e18 & 2.61e-9 & -2.63  & 2.34e-3 & -0.02 \\
 11 &  28673  &1e19 & 1.59e-8 &  2.36  & 7.12e-3 &  1.45 \\
 12 &  61441  &4e20 & 41      &  28    & 8e8     &  33   \\[0.5ex]
 \hline
\end{tabular}
\end{center}
\caption{Multilevel sparse collocation compared with sparse collocation using the Gaussian for Example \ref{EE 2D 1} when $C=3$. Max error evaluated at 64,000 Halton points in the whole domain.} \label{Table EE 2D 1 C=3 G MuSIK}
\end{table}

\begin{table}[!htb]
\begin{center}
\addtolength{\tabcolsep}{-1pt}
\begin{tabular}{||c|c|c|c|c||}
  \hline
\ Level \ & \ $E_{\rm extra} (C=2)$  \ & \ $\rho_{\rm extra} (C=2)$ \ & \ $E_{\rm extra} (C=3)$ \  & \ $\rho_{\rm extra} (C=3)$ \  \\
\hline \hline
 3    & 3.24e-2 & ---    & 5.66e-3  & ---   \\
 4    & 2.94e-3 & -2.87  & 2.55e-4  & -3.71 \\
 5    & 4.58e-4 & -2.26  & 4.35e-5  & -2.15 \\
 6    & 1.37e-4 & -1.49  & 3.33e-6  & -3.18 \\
 7    & 3.08e-5 & -1.87  & 3.64e-7  & -2.78 \\
 8    & 6.71e-6 & -1.93  & 5.04e-8  & -2.51 \\
 9    & 1.70e-6 & -1.76  & 6.73e-9  & -2.58 \\
 10   & 4.20e-7 & -1.81  & 8.67e-10 & -2.65 \\
 11   & 1.02e-7 & -1.84  & 3.47e-9  & 1.81  \\
 12   & 2.18e-8 & -2.02  & 2.00e-8  & 2.30  \\[0.5ex]
 \hline
\end{tabular}
\end{center}
\caption{Extrapolations from the multilevel sparse collocation using the Gaussian for Example \ref{EE 2D 1} with different constants: $C=2$ and $C=3$. Max error evaluated at 64,000 Halton points in the whole domain.} \label{Table EE 2D 1 Gauss MuSIK extra}
\end{table}

In Tables \ref{Table EE 2D 1 C=2 G MuSIK} and \ref{Table EE 2D 1 C=3 G MuSIK}, we use the Gaussian basis function in place of the MQ. Both tables also demonstrate the superiority of the multilevel method and that SIK-C using the Gaussian does not converge. When $C=3$, the condition number of MuSIK-C with the Gaussian reaches $1e19$ and $4e20$ at Levels $11$ and $12$, and the performance of MuSIK-C with the Gaussian breaks down. However, the condition number is $3e18$ at Level $10$ and the approximation is still improving. If we compare the condition number at which MQ approximation breaks down, around $1e14$, we conclude that the condition number alone is not a reliable indicator of the success of the method, and this depends on the basis function used. Similarly, in Table \ref{Table EE 2D 1 Gauss MuSIK extra}, we see better approximations using extrapolation, but the convergence rate is not increased.

\begin{figure}[htb!]
\centering
\begin{minipage}[b]{0.49\linewidth}
 \includegraphics[width=1\textwidth]{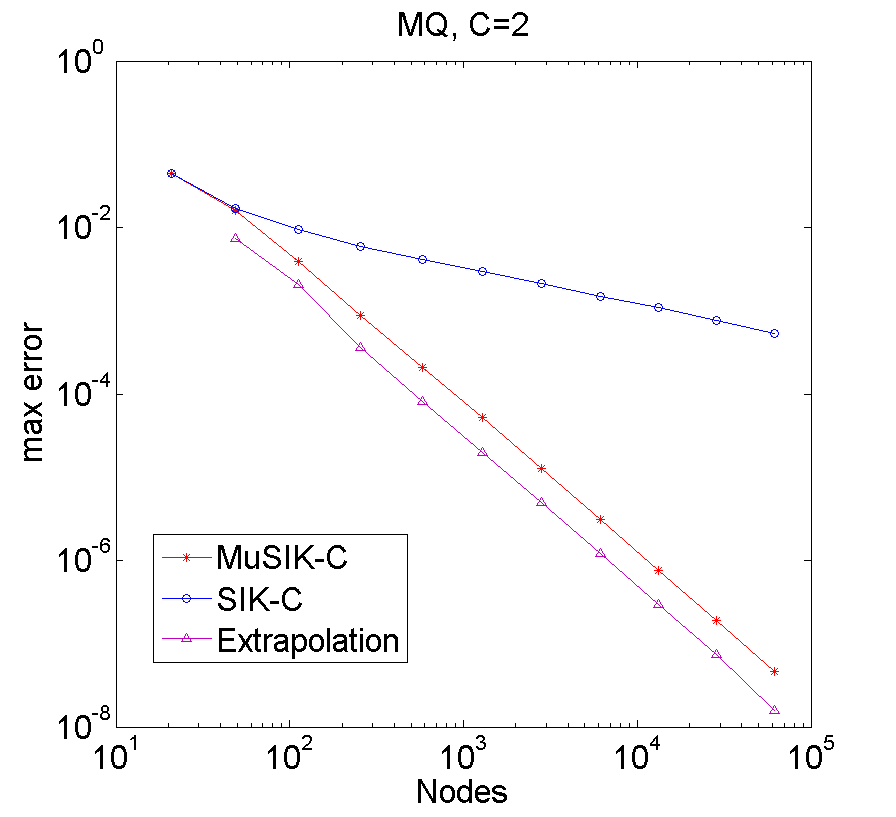}
\end{minipage}
\begin{minipage}[b]{0.49\linewidth}
 \includegraphics[width=1\textwidth]{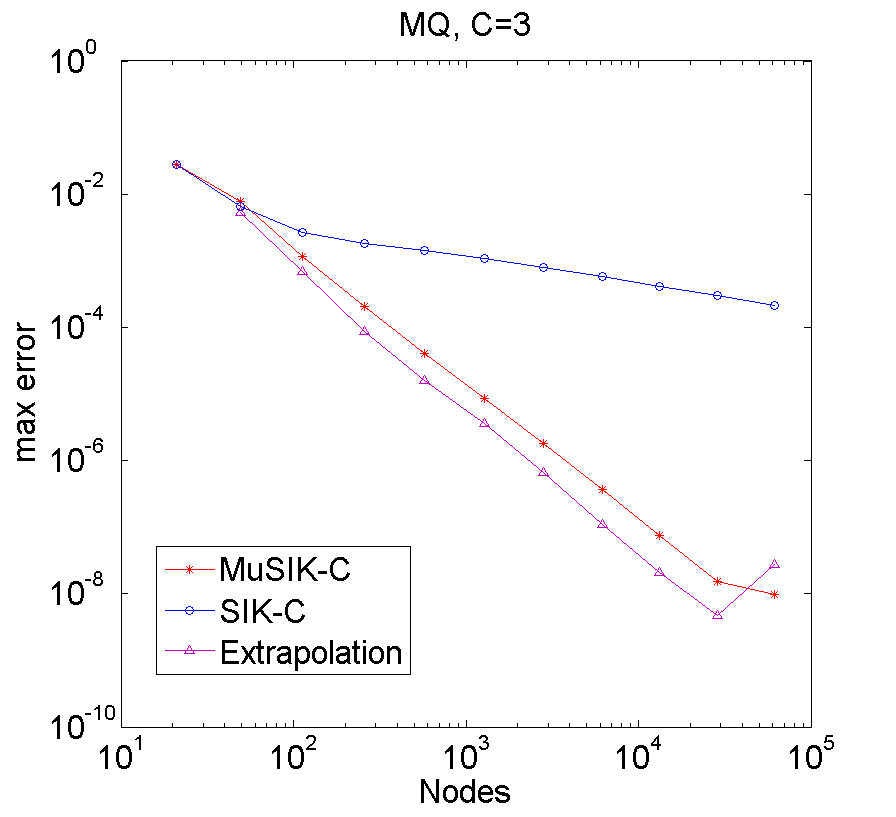}
\end{minipage}
\begin{minipage}[b]{0.49\linewidth}
 \includegraphics[width=1\textwidth]{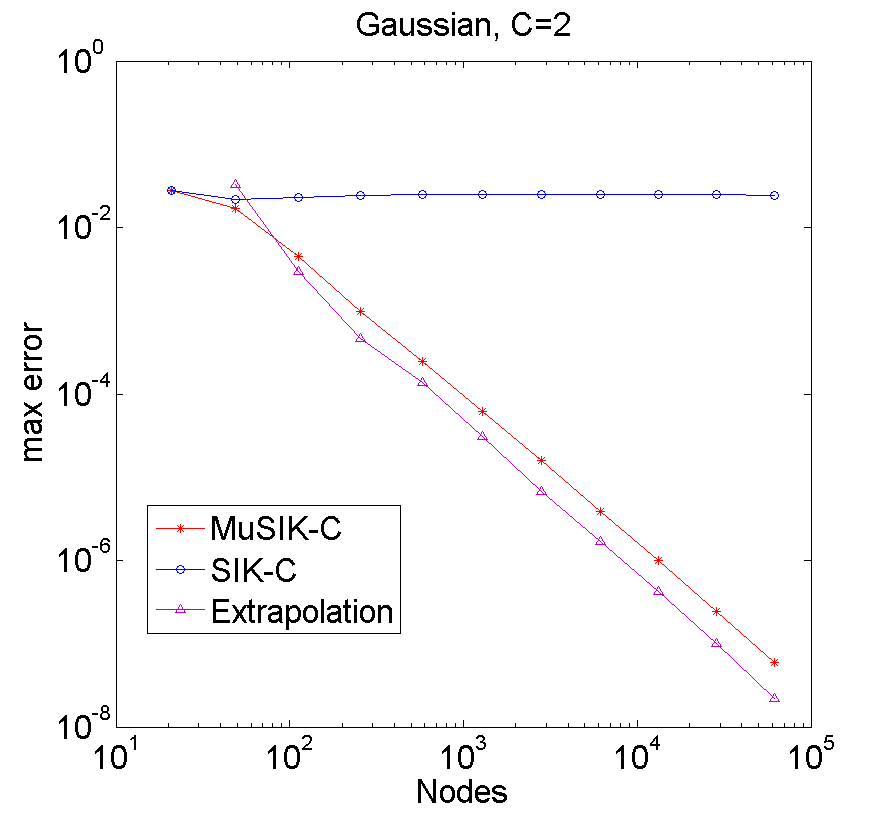}
\end{minipage}
\begin{minipage}[b]{0.49\linewidth}
 \includegraphics[width=1\textwidth]{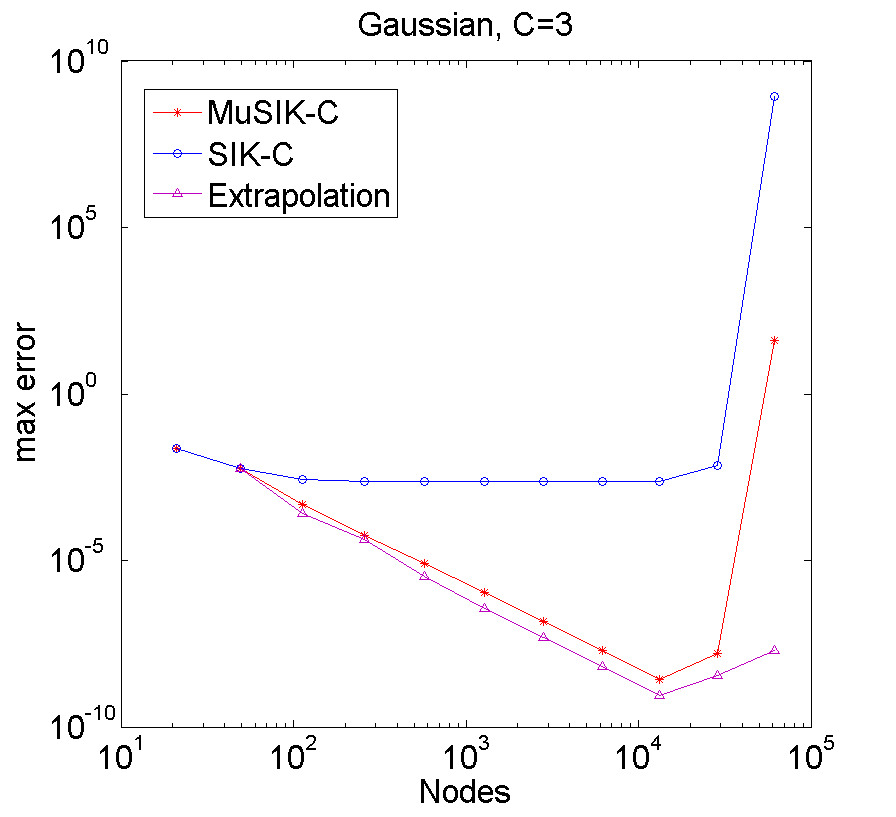}
\end{minipage}
\caption{The performance of SIK-C, MuSIK-C, and extrapolations with different basis functions and shape parameters for Example \ref{EE 2D 1}.} \label{Fig:EE 2D 1}
\end{figure}

In Figure \ref{Fig:EE 2D 1}, we observe that SIK-C  converges slowly using the MQ basis function, while for the Gaussian it does not converge at all. This phenomenon was also observed when approximating with Gaussian RBFs; see \cite{Usta}. We also observe the improvement from extrapolation. This figure illustrates nicely the advantages of using the multilevel method. In the remaining examples, we only show the solutions from MuSIK-C and extrapolations from MuSIK-C.

\clearpage

\begin{example} \label{EE 2D 2}
In this example, we solve the following two-dimensional problem on $\Omega=\left( 0,1 \right)^{2}$
\begin{equation}
\Delta u(\mathbf{x}) = - 2 \pi^{2} \sin(\pi x_{1}) \cos(\pi x_{2}) , \quad \mathbf{x} \in \Omega,
\end{equation}
with boundary conditions
\begin{equation}
u(\mathbf{x}) = \sin(\pi x_{1}) \cos(\pi x_{2}), \quad \mathbf{x} \in \partial \Omega.
\end{equation}
The exact solution is a two-dimensional tensor product function
\begin{equation}
u(\mathbf{x}) = \sin(\pi x_{1}) \cos(\pi x_{2}).
\end{equation}
\end{example}

\begin{table}[!htb]
\begin{center}
\addtolength{\tabcolsep}{-1pt}
\begin{tabular}{||c|c|c|c|c|c||}
  \hline
\ Level \ & \ Nodes \ & \ $E_{\rm MuSIK-C}$ \ & \ $\rho_{\rm MuSIK-C}$ \ & \ $E_{\rm Extra}$ \  & \ $\rho_{\rm Extra}$ \  \\
\hline \hline
 2  &  21     & 2.81e-2 & ---    & ---     & ---  \\
 3  &  49     & 7.53e-3 & -1.55  & 4.89e-3 & ---  \\
 4  &  113    & 2.07e-3 & -1.55  & 1.16e-3 & -1.73 \\
 5  &  257    & 5.33e-4 & -1.65  & 2.09e-4 & -2.08 \\
 6  &  577    & 1.33e-4 & -1.71  & 5.36e-5 & -1.68 \\
 7  &  1281   & 3.31e-5 & -1.75  & 1.33e-5 & -1.75 \\
 8  &  2817   & 8.14e-6 & -1.78  & 3.29e-6 & -1.77 \\
 9  &  6145   & 1.99e-6 & -1.80  & 8.04e-7 & -1.81 \\
 10 &  13313  & 4.89e-7 & -1.82  & 1.93e-7 & -1.84 \\
 11 &  28673  & 1.17e-7 & -1.86  & 4.90e-8 & -1.79 \\
 12 &  61441  & 2.83e-8 & -1.86  & 1.17e-8 & -1.88 \\[0.5ex]
 \hline
\end{tabular}
\end{center}
\caption{The performance of the multilevel sparse collocation method and corresponding extrapolations using the MQ for Example \ref{EE 2D 2} with $C=2$. Max error evaluated at 64,000 Halton points in the whole domain. $\rho_h \approx 2.0$.} \label{Table EE 2D 2 C=2 MQ MuSIK Extra}
\end{table}

\begin{table}[!htb]
\begin{center}
\addtolength{\tabcolsep}{-1pt}
\begin{tabular}{||c|c|c|c|c|c||}
  \hline
\ Level \ & \ Nodes \ & \ $E_{\rm MuSIK-C}$ \ & \ $\rho_{\rm MuSIK-C}$ \ & \ $E_{\rm Extra}$ \  & \ $\rho_{\rm Extra}$ \  \\
\hline \hline
 2  &  21     & 1.24e-2 & ---    & ---     & ---  \\
 3  &  49     & 4.02e-3 & -1.33  & 5.96e-3 & ---  \\
 4  &  113    & 8.54e-4 & -1.85  & 3.73e-4 & -3.32 \\
 5  &  257    & 1.89e-4 & -1.83  & 8.42e-5 & -1.81 \\
 6  &  577    & 4.12e-5 & -1.88  & 2.04e-5 & -1.75 \\
 7  &  1281   & 8.59e-6 & -1.97  & 3.84e-6 & -2.09 \\
 8  &  2817   & 1.74e-6 & -2.03  & 6.82e-7 & -2.19 \\
 9  &  6145   & 3.45e-7 & -2.07  & 1.22e-7 & -2.20 \\
 10 &  13313  & 6.69e-8 & -2.12  & 2.51e-8 & -2.05 \\
 11 &  28673  & 1.32e-8 & -2.11  & 4.68e-9 & -2.19 \\
 12 &  61441  & 6.72e-9 & -0.89  & 1.04e-8 &  1.05 \\[0.5ex]
 \hline
\end{tabular}
\end{center}
\caption{The performance of the multilevel sparse collocation method and corresponding extrapolations using the MQ for Example \ref{EE 2D 2} with $C=3$. Max error evaluated at 64,000 Halton points in the whole domain. $\rho_h \approx 2.3$.} \label{Table EE 2D 2 C=3 MQ MuSIK Extra}
\end{table}

\begin{table}[!htb]
\begin{center}
\addtolength{\tabcolsep}{-1pt}
\begin{tabular}{||c|c|c|c|c|c||}
  \hline
\ Level \ & \ Nodes \ & \ $E_{\rm MuSIK-C}$ \ & \ $\rho_{\rm MuSIK-C}$ \ & \ $E_{\rm Extra}$ \  & \ $\rho_{\rm Extra}$ \  \\
\hline \hline
 2  &  21     & 1.09e-2 & ---    & ---     & ---  \\
 3  &  49     & 5.01e-3 & -0.91  & 8.72e-3 & ---  \\
 4  &  113    & 1.70e-3 & -1.29  & 1.64e-3 & -2.00 \\
 5  &  257    & 5.50e-4 & -1.38  & 4.09e-4 & -1.69 \\
 6  &  577    & 1.62e-4 & -1.51  & 9.95e-5 & -1.75 \\
 7  &  1281   & 4.47e-5 & -1.61  & 2.35e-5 & -1.81 \\
 8  &  2817   & 1.18e-5 & -1.69  & 5.44e-6 & -1.86 \\
 9  &  6145   & 3.01e-6 & -1.75  & 1.24e-6 & -1.90 \\
 10 &  13313  & 7.55e-7 & -1.79  & 2.93e-7 & -1.87 \\
 11 &  28673  & 1.89e-7 & -1.80  & 7.21e-8 & -1.83 \\
 12 &  61441  & 4.69e-8 & -1.83  & 1.84e-8 & -1.79 \\[0.5ex]
 \hline
\end{tabular}
\end{center}
\caption{The performance of the multilevel sparse collocation method and corresponding extrapolations using the Gaussian for Example \ref{EE 2D 2} with $C=2$. Max error evaluated at 64,000 Halton points in the whole domain. $\rho_h \approx 2.0$} \label{Table EE 2D 2 C=2 G MuSIK Extra}
\end{table}

\begin{table}[!htb]
\begin{center}
\addtolength{\tabcolsep}{-1pt}
\begin{tabular}{||c|c|c|c|c|c||}
  \hline
\ Level \ & \ Nodes \ & \ $E_{\rm MuSIK-C}$ \ & \ $\rho_{\rm MuSIK-C}$ \ & \ $E_{\rm Extra}$ \  & \ $\rho_{\rm Extra}$ \  \\
\hline \hline
 2  &  21     & 1.20e-2  & ---    & ---      & ---  \\
 3  &  49     & 1.08e-3  & -2.85  & 1.39e-3  & ---  \\
 4  &  113    & 1.45e-4  & -2.40  & 7.64e-5  & -3.48 \\
 5  &  257    & 1.84e-5  & -2.51  & 7.14e-6  & -2.88 \\
 6  &  577    & 2.41e-6  & -2.51  & 8.42e-7  & -2.64 \\
 7  &  1281   & 3.25e-7  & -2.51  & 1.13e-7  & -2.52 \\
 8  &  2817   & 4.47e-8  & -2.52  & 1.43e-8  & -2.63 \\
 9  &  6145   & 5.98e-9  & -2.58  & 1.92e-9  & -2.57 \\
 10 &  13313  & 8.05e-10 & -2.59  & 2.86e-10 & -2.46 \\
 11 &  28673  & 3.64e-8  &  4.97  & 9.97e-10 &  1.63 \\
 12 &  61441  & 8.24     &  25.2  & 5.33e-8  &  5.22 \\[0.5ex]
 \hline
\end{tabular}
\end{center}
\caption{The performance of the multilevel sparse collocation method and corresponding extrapolations using the Gaussian for Example \ref{EE 2D 2} with $C=3$. Max error evaluated at 64,000 Halton points in the whole domain. $\rho_h \approx 2.9$.} \label{Table EE 2D 2 C=3 G MuSIK Extra}
\end{table}

\begin{figure}[!htb]
\begin{center}
   \includegraphics[width=6cm]{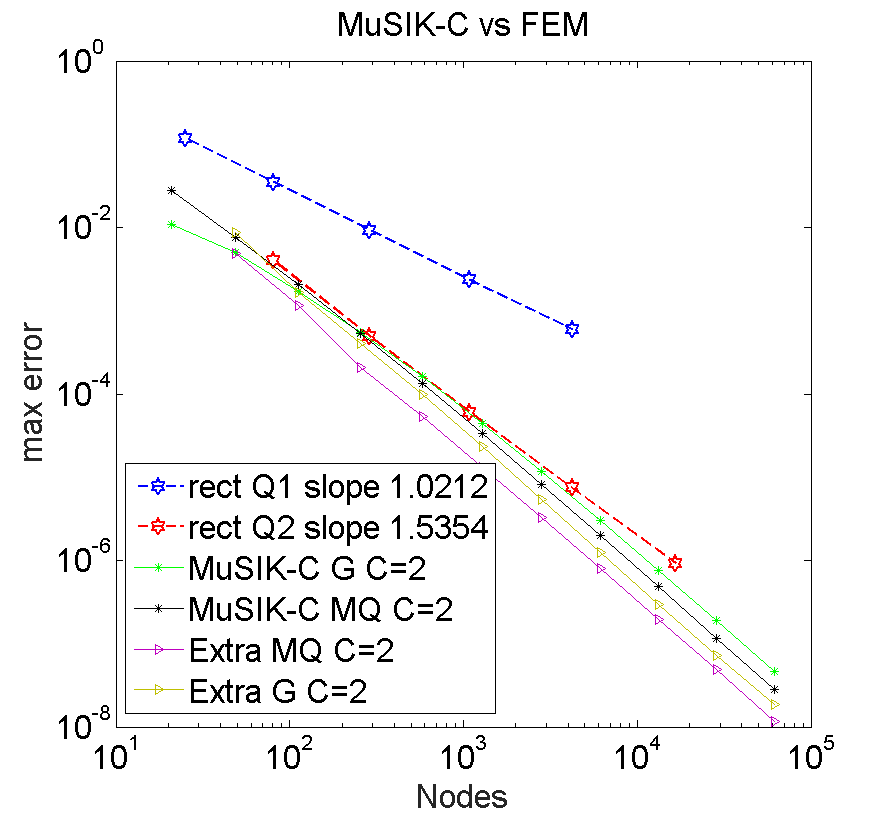}
   \includegraphics[width=6cm]{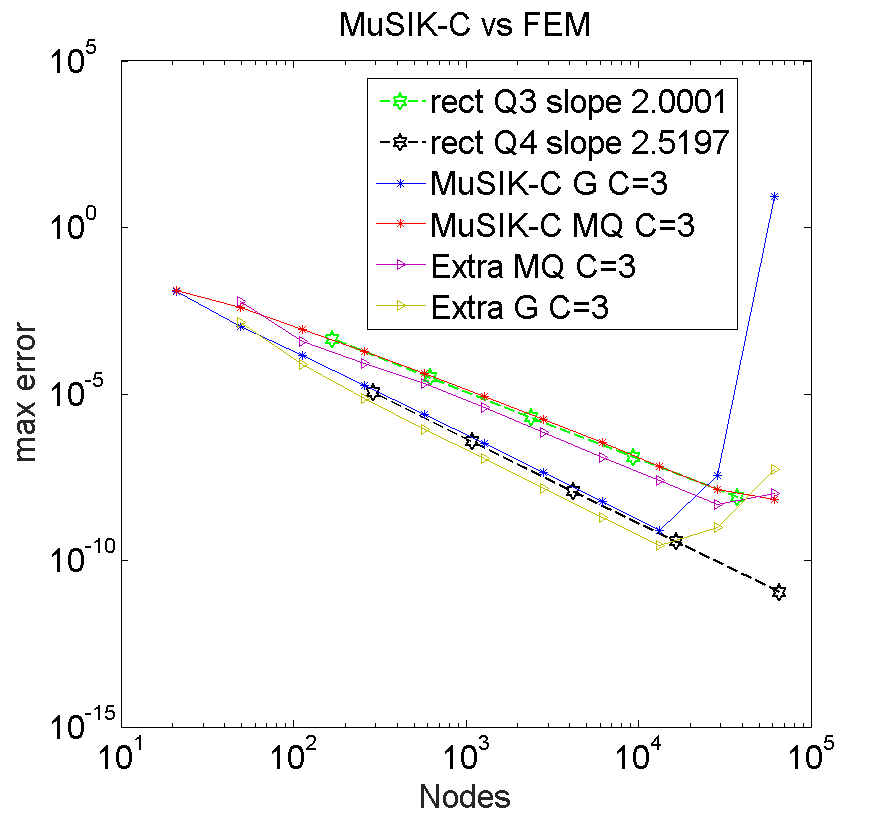}
\end{center}
\caption{The performance of MuSIK-C, extrapolation and FEM for Example \ref{EE 2D 2}.} \label{Fig:EE 2D 2}
\end{figure}

In Figure \ref{Fig:EE 2D 2}, the FEM is implemented by Z. Dong \cite{Dong}. The FEM is based on the Q-basis tensor product polynomials on a full grid. In Figure \ref{Fig:EE 2D 2} $Qp$ refers to a degree $p$ polynomial in each direction for the Q-basis method. We can see that the slope for the dashed line is almost growing as $\frac{p+1}{2}$. That means the convergence order is increasing with polynomial order $p$. Similarly, the convergence rate of MuSIK-C can be accelerated by increasing the shape parameter. In the left figure, the performance of MuSIK-C using both basis functions with $C=2$ is better than the performance of the FEM with $p=2$, for accuracy and convergence rate. In the right figure, with constant $C=3$, MuSIK-C using the MQ has similar performance with the FEM with $p=3$. Moreover, MuSIK-C using the Gaussian has similar performance with the FEM with $p=4$. However, MuSIK-C using both basis functions with $C=3$ breaks down at the last two levels because of the ill-condition. The convergence rate increases slowly before the ill-conditioning problem arises. This phenomenon demonstrates it is quite significant to reduce the condition number while utilising our MuSIK-C method. Extrapolation method always gives better approximations and does not improve the convergence rate.

\clearpage

\begin{example} \label{EE 3D 1}
In this example, we solve the following three-dimensional problem on $\Omega=\left( 0,1 \right)^{3}$
\begin{equation}
\Delta u(\mathbf{x}) = 0 , \quad \mathbf{x} \in \Omega,
\end{equation}
with boundary conditions
\begin{equation}
u(\mathbf{x}) = \sin(\pi x_{1}) \sin(\pi x_{2}) \frac{\sinh( \sqrt{2} \pi x_{3})}{\sinh(\sqrt{2} \pi)}, \quad \mathbf{x} \in \partial \Omega.
\end{equation}
The exact solution is the same as the boundary condition. 
\end{example}

\begin{figure}[htb!]
\begin{center}
   \includegraphics[width=6cm]{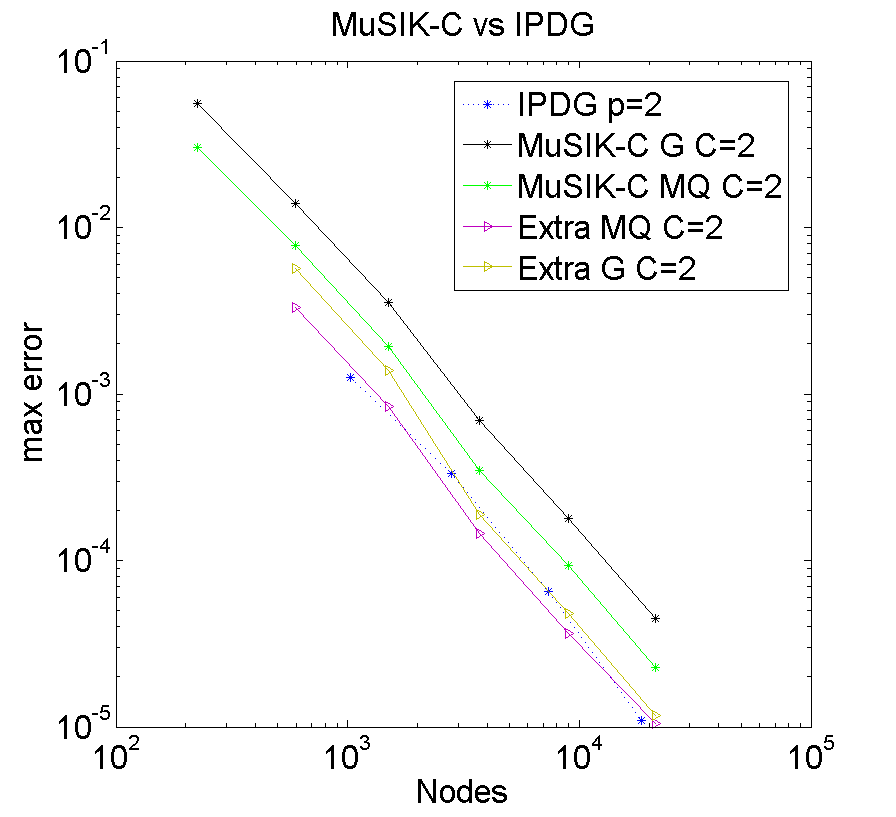}
   \includegraphics[width=6cm]{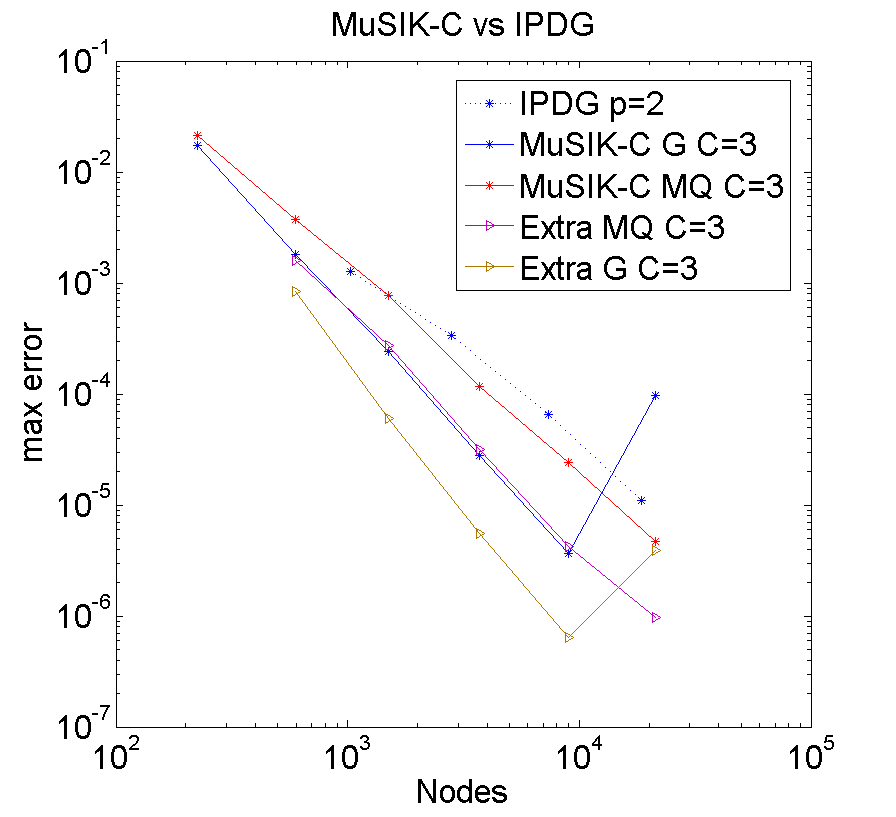}
\end{center}
\caption{Comparison between MuSIK-C and sparse grid IPDG for Example \ref{EE 3D 1}.} \label{Fig:EE 3D 1}
\end{figure}

In 2015, Wang et al. \cite{Wang1} developed an interior penalty discontinuous Galerkin (IPDG) method based on sparse grids to solve high-dimensional elliptic problems. Example \ref{EE 3D 1} and the corresponding numerical results for IPDG are all taken from \cite{Wang1}. Since IPDG is a sparse grid implementation the number of nodes used is directly comparable to ours. Figure \ref{Fig:EE 3D 1} shows the performance of the IPDG and our MuSIK-C for Examples \ref{EE 3D 1}.

\begin{table}[!htb]
\begin{center}
\addtolength{\tabcolsep}{-1pt}
\begin{tabular}{||c|c|c|c|c|c|c||}
  \hline
\ Level \ & \ Nodes \ & \ Cond \ & \  $E_{\rm MuSIK-C}$ \ & \ $\rho_{\rm MuSIK-C}$ \ & \ $E_{\rm Extra}$ \  & \ $\rho_{\rm Extra}$ \  \\
\hline \hline
 3 &  225    & 3e8  & 3.01e-2 & ---    & ---     & ---   \\
 4 &  593    & 2e9  & 7.83e-3 & -1.39  & 3.29e-3 & ---   \\
 5 &  1505   & 2e10 & 1.92e-3 & -1.51  & 8.42e-4 & -1.46 \\
 6 &  3713   & 1e11 & 3.49e-4 & -1.89  & 1.45e-4 & -1.95 \\
 7 &  8961   & 6e11 & 9.30e-5 & -1.50  & 3.66e-5 & -1.56 \\
 8 &  21249  & 3e12 & 2.29e-5 & -1.62  & 1.05e-5 & -1.44 \\ [0.5ex]
 \hline
\end{tabular}
\end{center}
\caption{The performance of the multilevel sparse collocation method and corresponding extrapolations using the MQ for Example \ref{EE 3D 1} with $C=2$. Max error evaluated at 120,000 Halton points in the whole domain. $\rho_h \approx 2.1$} \label{Table EE 3D 1 C=2 MQ MuSIK Extra}
\end{table}

\begin{table}[!htb]
\begin{center}
\addtolength{\tabcolsep}{-1pt}
\begin{tabular}{||c|c|c|c|c|c|c||}
  \hline
\ Level \ & \ Nodes \ & \ Cond \ & \  $E_{\rm MuSIK-C}$ \ & \ $\rho_{\rm MuSIK-C}$ \ & \ $E_{\rm Extra}$ \  & \ $\rho_{\rm Extra}$ \  \\
\hline \hline
 3 &  225    & 4e10 & 2.15e-2 & ---    & ---     & ---   \\
 4 &  593    & 7e11 & 3.74e-3 & -1.80  & 1.58e-3 & ---   \\
 5 &  1505   & 1e13 & 7.64e-4 & -1.70  & 2.73e-4 & -1.89 \\
 6 &  3713   & 1e14 & 1.17e-4 & -2.08  & 3.17e-5 & -2.39 \\
 7 &  8961   & 8e14 & 2.39e-5 & -1.80  & 4.22e-6 & -2.29 \\
 8 &  21249  & 6e15 & 4.67e-6 & -1.89  & 9.74e-7 & -1.70 \\ [0.5ex]
 \hline
\end{tabular}
\end{center}
\caption{The performance of the multilevel sparse collocation method and corresponding extrapolations using the MQ for Example \ref{EE 3D 1} with $C=3$. Max error evaluated at 120,000 Halton points in the whole domain. $\rho_h \approx 2.4$.} \label{Table EE 3D 1 C=3 MQ MuSIK Extra}
\end{table}

\begin{table}[!htb]
\begin{center}
\addtolength{\tabcolsep}{-1pt}
\begin{tabular}{||c|c|c|c|c|c|c||}
  \hline
\ Level \ & \ Nodes \ & \ Cond \ & \  $E_{\rm MuSIK-C}$ \ & \ $\rho_{\rm MuSIK-C}$ \ & \ $E_{\rm Extra}$ \  & \ $\rho_{\rm Extra}$ \  \\
\hline \hline
 3 &  225    & 1e7  & 5.55e-2 & ---    & ---     & ---   \\
 4 &  593    & 2e8  & 1.39e-2 & -1.43  & 5.65e-3 & ---   \\
 5 &  1505   & 1e9  & 3.52e-3 & -1.47  & 1.38e-3 & -1.52 \\
 6 &  3713   & 9e9  & 6.89e-4 & -1.81  & 1.89e-4 & -2.20 \\
 7 &  8961   & 6e10 & 1.78e-4 & -1.54  & 4.78e-5 & -1.56 \\
 8 &  21249  & 3e11 & 4.49e-5 & -1.60  & 1.17e-5 & -1.63 \\ [0.5ex]
 \hline
\end{tabular}
\end{center}
\caption{The performance of the multilevel sparse collocation method and corresponding extrapolations using the Gaussian for Example \ref{EE 3D 1} with $C=2$. Max error evaluated at 120,000 Halton points in the whole domain. $\rho_h \approx 2.0$.} \label{Table EE 3D 1 C=2 G MuSIK Extra}
\end{table}

\begin{table}[!htb]
\begin{center}
\addtolength{\tabcolsep}{-1pt}
\begin{tabular}{||c|c|c|c|c|c|c||}
  \hline
\ Level \ & \ Nodes \ & \ Cond \ & \  $E_{\rm MuSIK-C}$ \ & \ $\rho_{\rm MuSIK-C}$ \ & \ $E_{\rm Extra}$ \  & \ $\rho_{\rm Extra}$ \  \\
\hline \hline
 3 &  225    & 4e10 & 1.72e-2 & ---    & ---     & ---   \\
 4 &  593    & 1e13 & 1.79e-3 & -2.34  & 8.34e-4 & ---   \\
 5 &  1505   & 8e14 & 2.43e-4 & -2.14  & 6.01e-5 & -2.82 \\
 6 &  3713   & 3e16 & 2.78e-5 & -2.40  & 5.59e-6 & -2.63 \\
 7 &  8961   & 8e18 & 3.67e-6 & -2.30  & 6.36e-7 & -2.47 \\
 8 &  21249  & 5e20 & 9.74e-5 &  3.80  & 3.87e-6 &  2.09 \\ [0.5ex]
 \hline
\end{tabular}
\end{center}
\caption{The performance of the multilevel sparse collocation method and corresponding extrapolations using the Gaussian for Example \ref{EE 3D 1} with $C=3$. Max error evaluated at 120,000 Halton points in the whole domain. $\rho_h \approx 2.9$.} \label{Table EE 3D 1 C=3 G MuSIK Extra}
\end{table}

In the left figure of Figure \ref{Fig:EE 3D 1}, the MuSIK-method with $C=2$ (Gaussian and MQ) performs worse than the IPDG with polynomial order $p=2$. Extrapolation improves the accuracy so that we outperform IPDG with $p=2$. MuSIK-C with $C=3$ converges more quickly than IPDG while $p=2$. Again, extrapolation improves the performance of MuSIK-C. However, again we see that with the the higher value of $C$ we get ill-conditioning issues.

\begin{example} \label{EE New 3D}
In this example, we solve the following non tensor-product three-dimensional problem on $\Omega=\left( 0,1 \right)^{3}$
\begin{equation}
\Delta u(\mathbf{x}) = - \pi^{2} \sin \left( \pi \prod_{i=1}^{3} x_{i} \right) \left( \sum_{k=1}^{3} \prod_{j=1,j \neq k}^{3} x_{j} \right) , \quad \mathbf{x} \in \Omega,
\end{equation}
with boundary conditions
\begin{equation}
u(\mathbf{x}) = \sin \left( \pi \prod_{i=1}^{3} x_{i} \right), \quad \mathbf{x} \in \partial \Omega.
\end{equation}
The exact solution is a three-dimensional non-tensor product function
\begin{equation}
u(\mathbf{x}) = \sin \left( \pi \prod_{i=1}^{3} x_{i} \right).
\end{equation}
\end{example}

\begin{table}[!htb]
\begin{center}
\addtolength{\tabcolsep}{-1pt}
\begin{tabular}{||c|c|c|c|c|c|c||}
  \hline
\ Level \ & \ Nodes \ & \ Cond \ & \  $E_{\rm MuSIK-C}$ \ & \ $\rho_{\rm MuSIK-C}$ \ & \ $E_{\rm Extra}$ \  & \ $\rho_{\rm Extra}$ \  \\
\hline \hline
 3 &  225    & 3e8  & 2.35e-2 & ---    & ---     & ---   \\
 4 &  593    & 2e9  & 4.63e-3 & -1.68  & 4.63e-3 & ---   \\
 5 &  1505   & 2e10 & 9.66e-4 & -1.68  & 1.23e-3 & -1.42 \\
 6 &  3713   & 1e11 & 3.68e-4 & -1.07  & 8.32e-4 & -0.44 \\
 7 &  8961   & 6e11 & 1.83e-4 & -0.79  & 2.31e-4 & -1.46 \\
 8 &  21249  & 3e12 & 5.94e-5 & -1.30  & 3.48e-5 & -2.19 \\ [0.5ex]
 \hline
\end{tabular}
\end{center}
\caption{The performance of the multilevel sparse collocation method and corresponding extrapolations using the MQ for Example \ref{EE New 3D} with $C=2$. Max error evaluated at 120,000 Halton points in the whole domain. $\rho_h \approx 1.6$.} \label{Table EE New 3D C=2 MQ MuSIK Extra}
\end{table}

\begin{table}[!htb]
\begin{center}
\addtolength{\tabcolsep}{-1pt}
\begin{tabular}{||c|c|c|c|c|c|c||}
  \hline
\ Level \ & \ Nodes \ & \ Cond \ & \  $E_{\rm MuSIK-C}$ \ & \ $\rho_{\rm MuSIK-C}$ \ & \ $E_{\rm Extra}$ \  & \ $\rho_{\rm Extra}$ \  \\
\hline \hline
 3 &  225    & 4e10 & 1.17e-2 & ---    & ---     & ---   \\
 4 &  593    & 7e11 & 3.29e-3 & -1.31  & 5.16e-3 & ---   \\
 5 &  1505   & 1e13 & 8.99e-4 & -1.39  & 9.75e-4 & -1.79 \\
 6 &  3713   & 1e14 & 1.40e-4 & -2.06  & 1.35e-4 & -2.19 \\
 7 &  8961   & 8e14 & 2.72e-5 & -1.86  & 3.55e-5 & -1.52 \\
 8 &  21249  & 6e15 & 6.04e-6 & -1.74  & 6.85e-6 & -1.90 \\ [0.5ex]
 \hline
\end{tabular}
\end{center}
\caption{The performance of the multilevel sparse collocation method and corresponding extrapolations using the MQ for Example \ref{EE New 3D} with $C=3$. Max error evaluated at 120,000 Halton points in the whole domain. $\rho_h \approx 2.2$.} \label{Table EE New 3D C=3 MQ MuSIK Extra}
\end{table}

\begin{table}[!htb]
\begin{center}
\addtolength{\tabcolsep}{-1pt}
\begin{tabular}{||c|c|c|c|c|c|c||}
  \hline
\ Level \ & \ Nodes \ & \ Cond \ & \  $E_{\rm MuSIK-C}$ \ & \ $\rho_{\rm MuSIK-C}$ \ & \ $E_{\rm Extra}$ \  & \ $\rho_{\rm Extra}$ \  \\
\hline \hline
 3 &  225    & 1e7  & 2.67e-2 & ---    & ---     & ---   \\
 4 &  593    & 2e8  & 1.15e-2 & -0.88  & 2.07e-2 & ---   \\
 5 &  1505   & 1e9  & 2.92e-3 & -1.47  & 2.58e-3 & -2.23 \\
 6 &  3713   & 9e9  & 7.19e-4 & -1.55  & 4.42e-4 & -1.96 \\
 7 &  8961   & 6e10 & 1.91e-4 & -1.50  & 1.08e-4 & -1.59 \\
 8 &  21249  & 3e11 & 5.02e-5 & -1.55  & 2.91e-5 & -1.52 \\ [0.5ex]
 \hline
\end{tabular}
\end{center}
\caption{The performance of the multilevel sparse collocation method and corresponding extrapolations using the Gaussian for Example \ref{EE New 3D} with $C=2$. Max error evaluated at 120,000 Halton points in the whole domain. $\rho_h \approx 1.9$.} \label{Table EE New 3D C=2 G MuSIK Extra}
\end{table}

\begin{table}[!htb]
\begin{center}
\addtolength{\tabcolsep}{-1pt}
\begin{tabular}{||c|c|c|c|c|c|c||}
  \hline
\ Level \ & \ Nodes \ & \ Cond \ & \  $E_{\rm MuSIK-C}$ \ & \ $\rho_{\rm MuSIK-C}$ \ & \ $E_{\rm Extra}$ \  & \ $\rho_{\rm Extra}$ \  \\
\hline \hline
 3 &  225    & 4e10 & 7.91e-3 & ---    & ---     & ---   \\
 4 &  593    & 1e13 & 1.65e-3 & -1.62  & 4.16e-3 & ---   \\
 5 &  1505   & 8e14 & 5.07e-4 & -1.26  & 6.66e-4 & -1.97 \\
 6 &  3713   & 3e16 & 6.95e-5 & -2.20  & 3.91e-5 & -3.14 \\
 7 &  8961   & 8e18 & 5.78e-6 & -2.82  & 2.60e-6 & -3.07 \\
 8 &  21249  & 5e20 & 2.36e-6 & -1.04  & 4.13e-6 &  0.54 \\ [0.5ex]
 \hline
\end{tabular}
\end{center}
\caption{The performance of the multilevel sparse collocation method and corresponding extrapolations using the Gaussian for Example \ref{EE New 3D} with $C=3$. Max error evaluated at 120,000 Halton points in the whole domain. $\rho_h \approx 3.6$.} \label{Table EE New 3D C=3 G MuSIK Extra}
\end{table}

\begin{figure}[htb!]
\begin{center}
   \includegraphics[width=6cm]{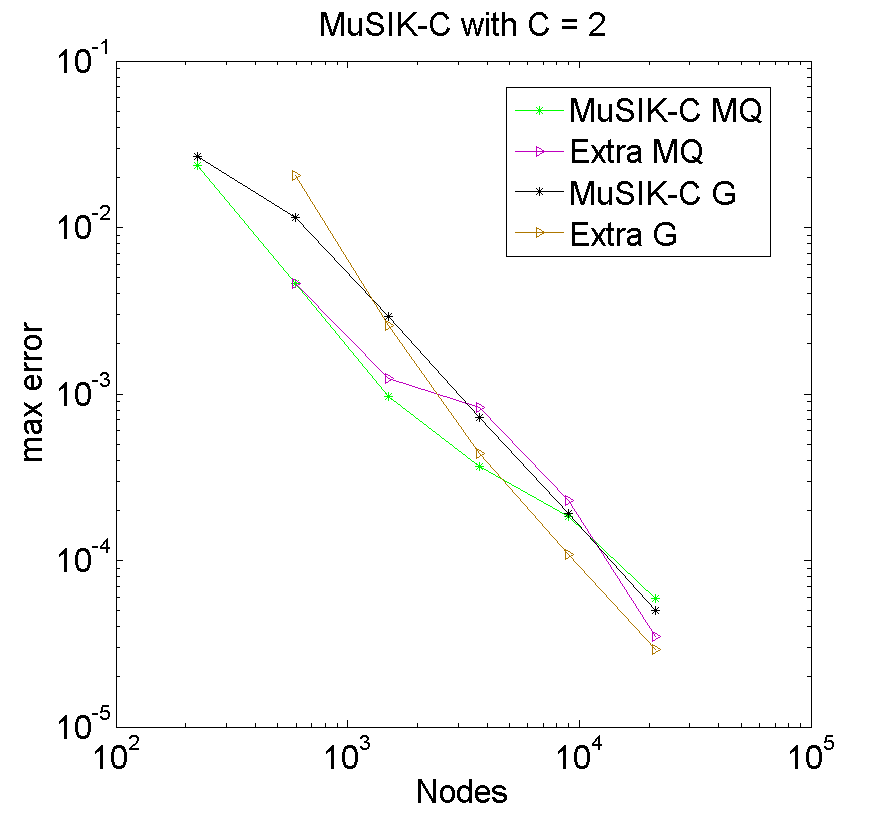}
   \includegraphics[width=6cm]{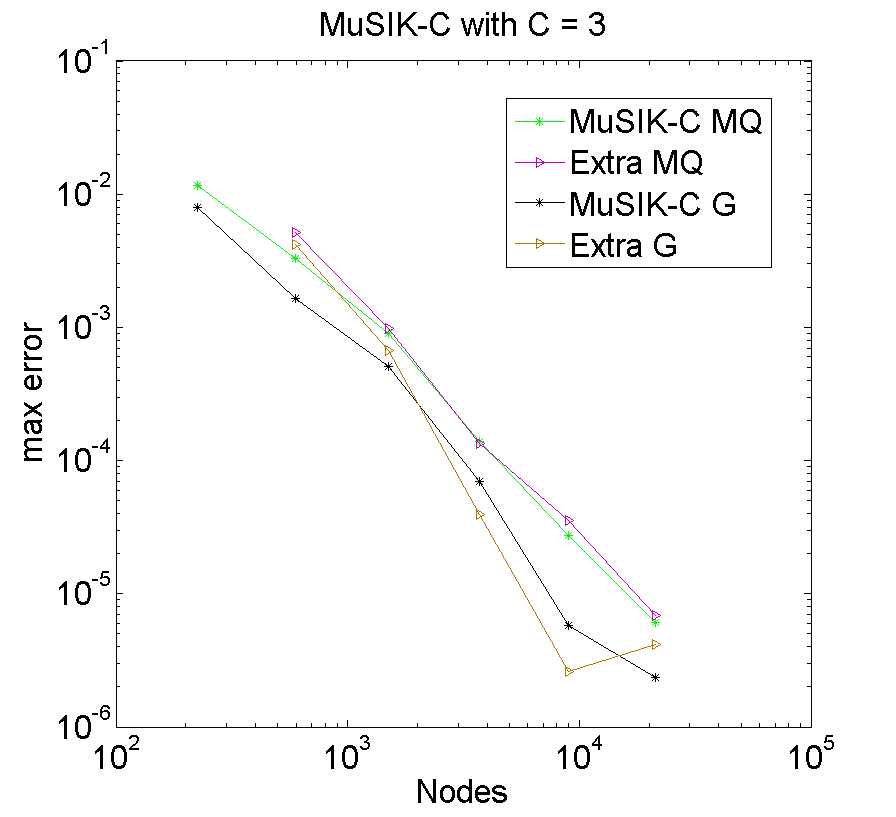}
\end{center}
\caption{The performance of MuSIK-C and extrapolation for Example \ref{EE New 3D}.} \label{Fig:EE New 3D}
\end{figure}

If we compare Tables~\ref{Table EE New 3D C=2 MQ MuSIK Extra}~to~\ref{Table EE New 3D C=3 G MuSIK Extra} to Tables~\ref{Table EE 3D 1 C=2 MQ MuSIK Extra}~to~\ref{Table EE 3D 1 C=3 G MuSIK Extra}, we see that the non tensor product convergence rates are poorer for the range of examples we could compute than the tensor product case.

\begin{example} \label{EE 4D 1}
In this example, we solve a four-dimensional problem on $\Omega=\left(0,1 \right)^{4}$
\begin{equation}
\Delta u(\mathbf{x}) = 0 , \quad \mathbf{x} \in \Omega,
\end{equation}
with boundary conditions
\begin{equation}
u(\mathbf{x}) = \sin(\pi x_{1}) \sin(\pi x_{2}) \sin(\pi x_{3}) \frac{\sinh( \sqrt{3} \pi x_{4})}{\sinh(\sqrt{3} \pi)}, \quad \mathbf{x} \in \partial \Omega,
\end{equation}
with the boundary condition as the exact solution.
\end{example}

\begin{table}[!htb]
\begin{center}
\begin{tabular}{||c|c|c|c|c|c|c||}
  \hline
\ Level \ & \ Nodes \ & \ Cond \ & \ $E_{\rm MuSIK-C}$ \ & \ $\rho_{\rm MuSIK-C}$ \ & \ $E_{\rm extra}$ \  & \ $\rho_{\rm extra}$ \  \\
\hline \hline
 4 &  2769   & 4e11 & 1.62e-2 & ---    & ---     & ---   \\
 5 &  7681   & 3e12 & 3.93e-3 & -1.39  & 1.74e-3 & ---   \\
 6 &  20481  & 3e13 & 9.62e-4 & -1.44  & 4.13e-4 & -1.47 \\
 7 &  52993  & 2e14 & 2.42e-4 & -1.45  & 1.28e-4 & -1.23 \\
 8 &  133889 & 1e15 & 5.91e-5 & -1.52  & 3.43e-5 & -1.42 \\ [0.5ex]
 \hline
\end{tabular}
\end{center}
\caption{The performance of the multilevel sparse collocation method and corresponding extrapolations using the MQ for Example \ref{EE 4D 1} with $C=2$. Max error evaluated at 240,000 Halton points in the whole domain. $\rho_h \approx 2.0$.} \label{Table EE 4D 1 C=2 MuSIK MQ}
\end{table}

\begin{table}[!htb]
\begin{center}
\begin{tabular}{||c|c|c|c|c|c|c||}
  \hline
\ Level \ & \ Nodes \ & \ Cond \ & \ $E_{\rm MuSIK-C}$ \ & \ $\rho_{\rm MuSIK-C}$ \ & \ $E_{\rm extra}$ \  & \ $\rho_{\rm extra}$ \  \\
\hline \hline
 4 &  2769   & 4e9  & 4.26e-2 & ---    & ---     & ---   \\
 5 &  7681   & 4e10 & 1.02e-2 & -1.40  & 4.01e-3 & ---   \\
 6 &  20481  & 3e11 & 2.52e-3 & -1.43  & 1.01e-3 & -1.40 \\
 7 &  52993  & 2e12 & 6.40e-4 & -1.44  & 2.44e-4 & -1.50 \\
 8 &  133889 & 2e13 & 1.58e-4 & -1.51  & 6.53e-5 & -1.42 \\ [0.5ex]
 \hline
\end{tabular}
\end{center}
\caption{The performance of the multilevel sparse collocation method and corresponding extrapolations using the Gaussian for Example \ref{EE 4D 1} with $C=2$. Max error evaluated at 240,000 Halton points in the whole domain. $\rho_h \approx 2.0$.} \label{Table EE 4D 1 C=2 MuSIK G}
\end{table}

\begin{figure}[htb!]
\centering
 \includegraphics[width=10cm]{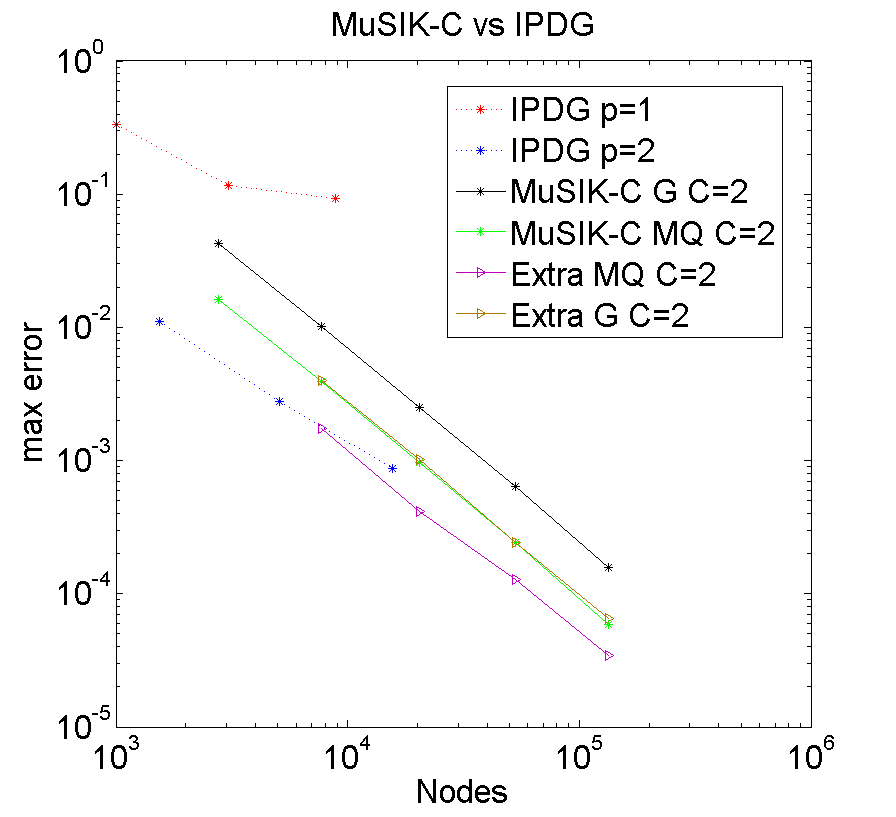}
\caption{Comparison between MuSIK-C and sparse grid IPDG for Example \ref{EE 4D 1}.} \label{Fig:EE 4D 1}
\end{figure}

We compare our results again with those from \cite{Wang1}. In Figure \ref{Fig:EE 4D 1}, we see that MuSIK-C with $C=2$ has a better convergence rate than IPDG with $p=2$ (note the steeper slopes). This suggests that MuSIK-C might improve relatively as dimension increases, but this assertion requires further experimentation. Again, extrapolation improves the approximation accuracy.

\begin{example} \label{EE New 4D}
In this example, we solve the following four-dimensional non tensor-product problem on $\Omega=\left(0,1 \right)^{4}$
\begin{equation}
\Delta u(\mathbf{x}) = - \pi^{2} \sin \left( \pi \prod_{j=1}^{4} x_{j} \right) \left( \sum_{k=1}^{4} \prod_{i=1,i \neq k}^{4} x_{i}  \right) , \quad \mathbf{x} \in \Omega,
\end{equation}
with boundary conditions
\begin{equation}
u(\mathbf{x}) = \sin \left( \pi \prod_{i=1}^{4} x_{i} \right), \quad \mathbf{x} \in \partial \Omega.
\end{equation}
\end{example}

\begin{table}[!htb]
\begin{center}
\begin{tabular}{||c|c|c|c|c|c|c||}
  \hline
\ Level \ & \ Nodes \ & \ Cond \ & \ $E_{\rm MuSIK-C}$ \ & \ $\rho_{\rm MuSIK-C}$ \ & \ $E_{\rm extra}$ \  & \ $\rho_{\rm extra}$ \  \\
\hline \hline
 4 &  2769   & 4e11 & 1.03e-2 & ---    & ---     & ---   \\
 5 &  7681   & 3e12 & 5.24e-3 & -0.66  & 1.59e-2 & ---   \\
 6 &  20481  & 3e13 & 1.86e-3 & -1.06  & 1.27e-3 & -2.58 \\
 7 &  52993  & 2e14 & 5.06e-4 & -1.37  & 2.27e-4 & -1.81 \\
 8 &  133889 & 1e15 & 1.37e-4 & -1.41  & 4.76e-5 & -1.68 \\ [0.5ex]
 \hline
\end{tabular}
\end{center}
\caption{The performance of the multilevel sparse collocation method and corresponding extrapolations using the MQ for Example \ref{EE New 4D} with $C=2$. Max error evaluated at 240,000 Halton points in the whole domain. $\rho_h \approx 1.8$.} \label{Table EE New 4D C=2 MuSIK MQ}
\end{table}

\begin{table}[!htb]
\begin{center}
\begin{tabular}{||c|c|c|c|c|c|c||}
  \hline
\ Level \ & \ Nodes \ & \ Cond \ & \ $E_{\rm MuSIK-C}$ \ & \ $\rho_{\rm MuSIK-C}$ \ & \ $E_{\rm extra}$ \  & \ $\rho_{\rm extra}$ \  \\
\hline \hline
 4 &  2769   & 4e9  & 2.84e-2 & ---    & ---     & ---   \\
 5 &  7681   & 4e10 & 1.16e-2 & -0.88  & 3.65e-2 & ---   \\
 6 &  20481  & 3e11 & 2.20e-3 & -1.70  & 1.84e-3 & -3.05 \\
 7 &  52993  & 2e12 & 6.98e-4 & -1.21  & 4.77e-4 & -1.42 \\
 8 &  133889 & 2e13 & 2.03e-4 & -1.33  & 9.86e-5 & -1.70 \\ [0.5ex]
 \hline
\end{tabular}
\end{center}
\caption{The performance of the multilevel sparse collocation method and corresponding extrapolations using the Gaussian for Example \ref{EE New 4D} with $C=2$. Max error evaluated at 240,000 Halton points in the whole domain. $\rho_h \approx 1.8$.} \label{Table EE New 4D C=2 MuSIK G}
\end{table}

\begin{figure}[htb!]
\centering
 \includegraphics[width=10cm]{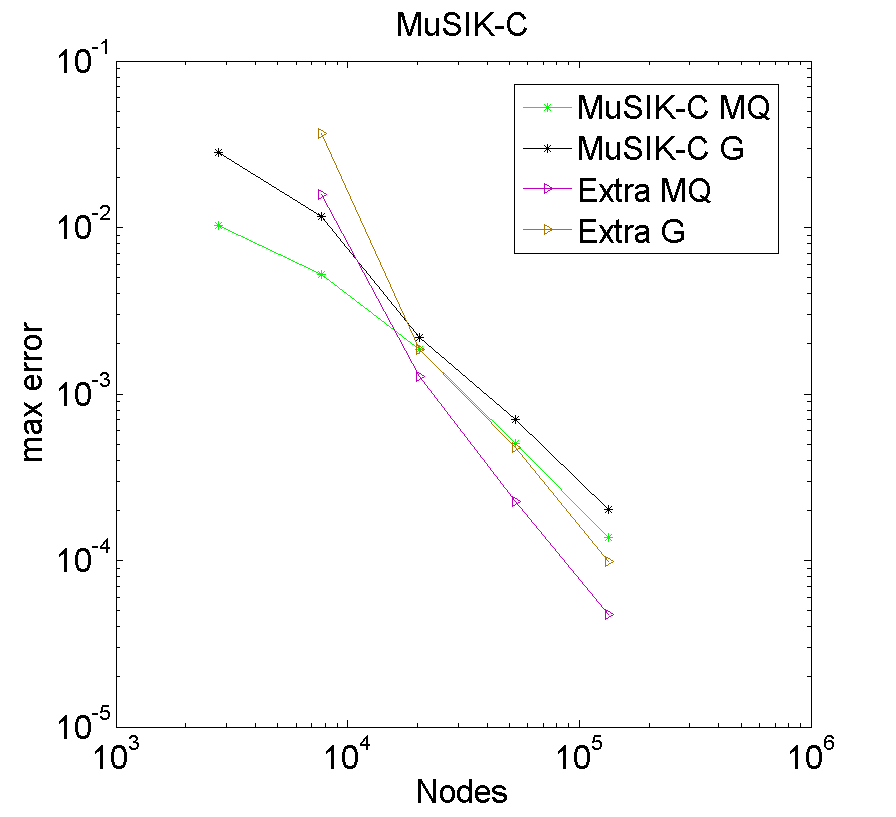}
\caption{The performance of MuSIK-C and extrapolation using $C = 2$ for Example \ref{EE New 4D}.} \label{Fig:EE New 4D}
\end{figure}

Comparison of Tables~\ref{Table EE New 4D C=2 MuSIK MQ} and \ref{Table EE New 4D C=2 MuSIK G} with Tables~\ref{Table EE 4D 1 C=2 MuSIK MQ} and \ref{Table EE 4D 1  C=2 MuSIK G} show that we have better convergence rates in for the tensor product case than for non tensor product cases in this instance.

\subsection{Parabolic examples}

Langer et al. \cite{Langer1} presented the new stable space-time Isogeometric Analysis (IgA) method in 2016. Isogeometric analysis is a collection of methods that use splines, or some of their extensions such as NURBS (non-uniform rational B-splines) and T-splines, as functions to build approximation spaces which are then used to solve partial differential equations numerically. As the authors just presented $L_{2}$ errors in \cite{Langer1}, we change to RMS error here to compare. Let us define the error and rate as
\begin{eqnarray*}
E^{n,d}_{\rm RMS} &=& \sqrt{ \frac{1}{N_{\mathbf{T}}} \sum^{N_{\mathbf{T}}}_{i=1} \left( u(\bx_{i}) - \hat{u}_{\rm ML}^{n,d}(\bx_{i}) \right)^{2} }, \quad \bx_{i} \in \mathbf{T}, \\
\rho_{\rm RMS}^{n+1,d} &=& \frac{ \log \left( E_{\rm RMS}^{n+1,d} \right) - \log \left( E_{\rm RMS}^{n,d} \right) }{\log \left( {\rm Nodes}^{n+1,d} \right) - \log \left( {\rm Nodes}^{n,d} \right)}.
\end{eqnarray*}

In these time-dependent examples we just apply an initial condition, the final boundary is left open. There is no computational issue for the method in doing this, and as we see below, the results are good.

\begin{example} \label{PE 4D 1}
In this example, we solve the following three-dimensional spatial problem on $\Omega_{t} = \Omega \times t  = \left[0,1 \right]^{3} \times \left[0,1 \right]$
\begin{equation}
u_{t} - \Delta u = \pi \sin (\pi x_{1}) \sin (\pi x_{2}) \sin (\pi x_{3}) \left( \cos (\pi t) + 3\pi \sin (\pi t) \right),\quad \bx \in \Omega, \quad t \in \left(0,1 \right],
\end{equation}
with boundary and initial conditions
\begin{eqnarray}
u(t,\mathbf{x}) &=& \sin(\pi t) \sin(\pi x_{1}) \sin(\pi x_{2}) \sin(\pi x_{3}),\quad \bx \in \partial \Omega,\quad t \in \left(0,1 \right],\\
u(0,\mathbf{x}) &=& 0,\quad \bx \in \Omega.
\end{eqnarray}
\end{example}

\begin{table}[!htb]
\begin{center}
\begin{tabular}{||c|c|c|c|c|c|c||}
  \hline
\ Level \ & \ Nodes \ & \ Cond \ & \ $E_{\rm RMS}$ \ & \ $\rho_{\rm RMS}$ \ & \ $E_{\rm extra}$ \  & \ $\rho_{\rm extra}$ \  \\
\hline \hline
 4 &  2769   & 3e11 & 1.80e-3 & ---    & ---     & ---   \\
 5 &  7681   & 3e12 & 3.71e-4 & -1.55  & 1.60e-4 & ---   \\
 6 &  20481  & 2e13 & 8.91e-5 & -1.45  & 2.12e-5 & -2.06 \\
 7 &  52993  & 2e14 & 2.18e-5 & -1.48  & 3.75e-6 & -1.82 \\
 8 &  133889 & 1e15 & 5.48e-6 & -1.49  & 1.19e-6 & -1.24 \\ [0.5ex]
 \hline
\end{tabular}
\end{center}
\caption{The performance of the multilevel sparse collocation method and corresponding extrapolation using the MQ for Example \ref{PE 4D 1} with $C=2$. Max error evaluated at 240,000 Halton points in the whole domain. $\rho_h \approx 2.0$} \label{Table PE 4D 1 C=2 MuSIK MQ}
\end{table}

\begin{table}[!htb]
\begin{center}
\begin{tabular}{||c|c|c|c|c|c|c||}
  \hline
\ Level \ & \ Nodes \ & \ Cond \ & \ $E_{\rm RMS}$ \ & \ $\rho_{\rm RMS}$ \ & \ $E_{\rm extra}$ \  & \ $\rho_{\rm extra}$ \  \\
\hline \hline
 4 &  2769   & 5e9  & 3.19e-3 & ---    & ---     & ---   \\
 5 &  7681   & 4e10 & 6.98e-4 & -1.49  & 2.31e-4 & ---   \\
 6 &  20481  & 3e11 & 1.70e-4 & -1.44  & 2.73e-5 & -2.18 \\
 7 &  52993  & 2e12 & 4.23e-5 & -1.46  & 5.26e-6 & -1.73 \\
 8 &  133889 & 3e13 & 1.06e-5 & -1.50  & 1.06e-6 & -1.73 \\ [0.5ex]
 \hline
\end{tabular}
\end{center}
\caption{The performance of the multilevel sparse collocation method and corresponding extrapolations using the Gaussian for Example \ref{PE 4D 1} with $C=2$. Max error evaluated at 240,000 Halton points in the whole domain. $\rho_h \approx 2.0$.} \label{Table PE 4D 1 C=2 MuSIK G}
\end{table}

\begin{figure}[htb!]
\centering
 \includegraphics[width=10cm]{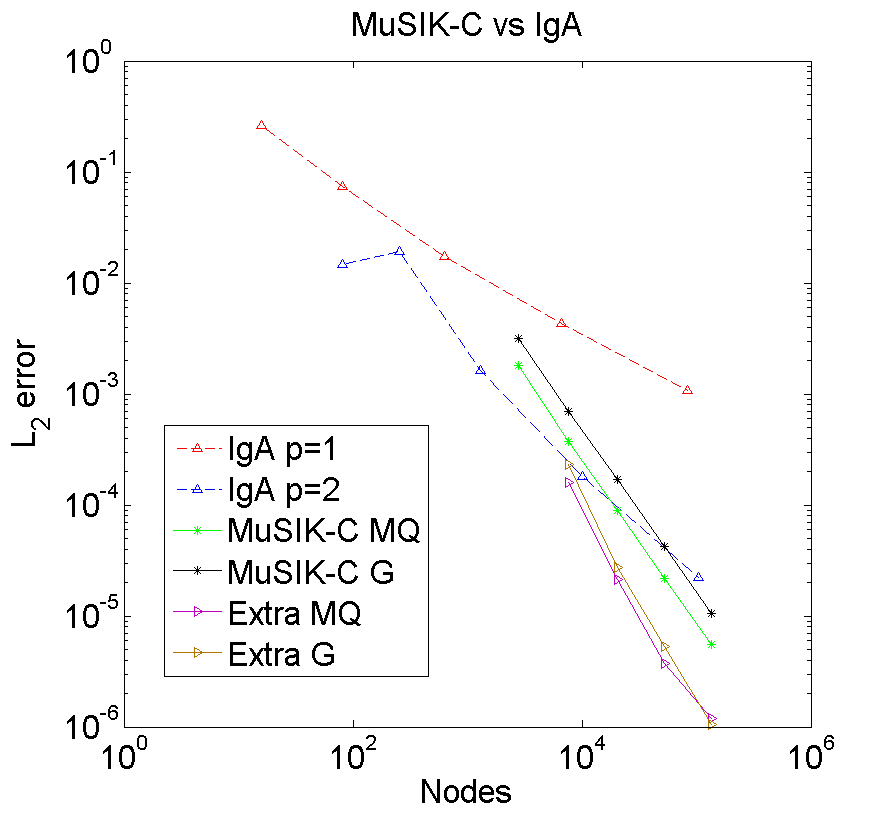}
\caption{Comparison between MuSIK-C and IgA for Example \ref{PE 4D 1}.} \label{Fig:PE 4D 1}
\end{figure}

In Tables \ref{Table PE 4D 1 C=2 MuSIK MQ} and \ref{Table PE 4D 1 C=2 MuSIK G}, we see that the condition numbers grow quickly for both basis functions. In these two tables, the focus is on multilevel sparse collocation using MQ and Gaussian with constant $C=2$. The convergence rate here also appears to grow slowly, which is indicative of faster than polynomial convergence rates. In Figure \ref{Fig:PE 4D 1}, the numerical results of IgA are taken from \cite{Langer1}, so we also present RMS errors here for comparison purposes. MuSIK-C with $C=2$ appears to have faster convergence than IgA with $p=1$ and $p=2$ on the hypeprcube. Since MuSIK-C is not applicable on an irregular domain, we cannot compare with IgA on such regions.

\begin{example} \label{PE 4D 2}
In this example, we solve the following three-dimensional spatial problem on $\Omega_{t} = \Omega \times t  = \left[0,1 \right]^{3} \times \left[0,1 \right]$
\begin{equation}
u_{t} - \Delta u = e^{10(t-1)}\sin(\pi x_{1} x_{2} x_{3}) \left( 10 + \pi^{2} \left( x_{2}^{2}x_{3}^{2} + x_{1}^{2}x_{3}^{2} + x_{1}^{2}x_{2}^{2} \right) \right),\quad \bx \in \Omega, \quad t \in \left(0,1 \right],
\end{equation}
with non tensor-product boundary and initial conditions
\begin{eqnarray}
u(t,\mathbf{x}) &=& e^{10(t-1)}\sin(\pi x_{1} x_{2} x_{3}),\quad \bx \in \partial \Omega,\quad t \in \left(0,1 \right],\\
u(0,\mathbf{x}) &=& e^{-10}\sin(\pi x_{1} x_{2} x_{3}),\quad \bx \in \Omega.
\end{eqnarray}
\end{example}

\begin{table}[!htb]
\begin{center}
\begin{tabular}{||c|c|c|c|c|c|c||}
  \hline
\ Level \ & \ Nodes \ & \ Cond \ & \ $E_{\rm MuSIK-C}$ \ & \ $\rho_{\rm MuSIK-C}$ \ & \ $E_{\rm extra}$ \  & \ $\rho_{\rm extra}$ \  \\
\hline \hline
 4 &  2769   & 3e11 & 6.44e-2 & ---    & ---     & ---   \\
 5 &  7681   & 3e12 & 2.70e-2 & -0.85  & 2.08e-2 & ---   \\
 6 &  20481  & 2e13 & 9.54e-3 & -1.06  & 5.37e-3 & -1.38 \\
 7 &  52993  & 2e14 & 3.30e-3 & -1.12  & 1.45e-3 & -1.38 \\
 8 &  133889 & 1e15 & 1.06e-3 & -1.22  & 3.99e-4 & -1.39 \\ [0.5ex]
 \hline
\end{tabular}
\end{center}
\caption{The performance of the multilevel sparse collocation method and corresponding extrapolations using the MQ for Example \ref{PE 4D 2} with $C=2$. Max error evaluated at 240,000 Halton points in the whole domain. $\rho_h \approx 1.6$} \label{Table PE 4D 2 C=2 MuSIK MQ}
\end{table}

\begin{table}[!htb]
\begin{center}
\begin{tabular}{||c|c|c|c|c|c|c||}
  \hline
\ Level \ & \ Nodes \ & \ Cond \ & \ $E_{\rm MuSIK-C}$ \ & \ $\rho_{\rm MuSIK-C}$ \ & \ $E_{\rm extra}$ \  & \ $\rho_{\rm extra}$ \  \\
\hline \hline
 4 &  2769   & 5e9  & 9.18e-2 & ---    & ---     & ---   \\
 5 &  7681   & 4e10 & 4.33e-2 & -0.74  & 4.07e-2 & ---   \\
 6 &  20481  & 3e11 & 1.30e-2 & -1.23  & 6.36e-3 & -1.89 \\
 7 &  52993  & 2e12 & 4.01e-3 & -1.24  & 1.66e-3 & -1.41 \\
 8 &  133889 & 3e13 & 1.18e-3 & -1.32  & 3.87e-4 & -1.57 \\ [0.5ex]
 \hline
\end{tabular}
\end{center}
\caption{The performance of the multilevel sparse collocation method and corresponding extrapolations using the Gaussian for Example \ref{PE 4D 2} with $C=2$. Max error evaluated at 240,000 Halton points in the whole domain. $\rho_h \approx 1.76$.} \label{Table PE 4D 2 C=2 MuSIK G}
\end{table}

\begin{figure}[htb!]
\centering
 \includegraphics[width=10cm]{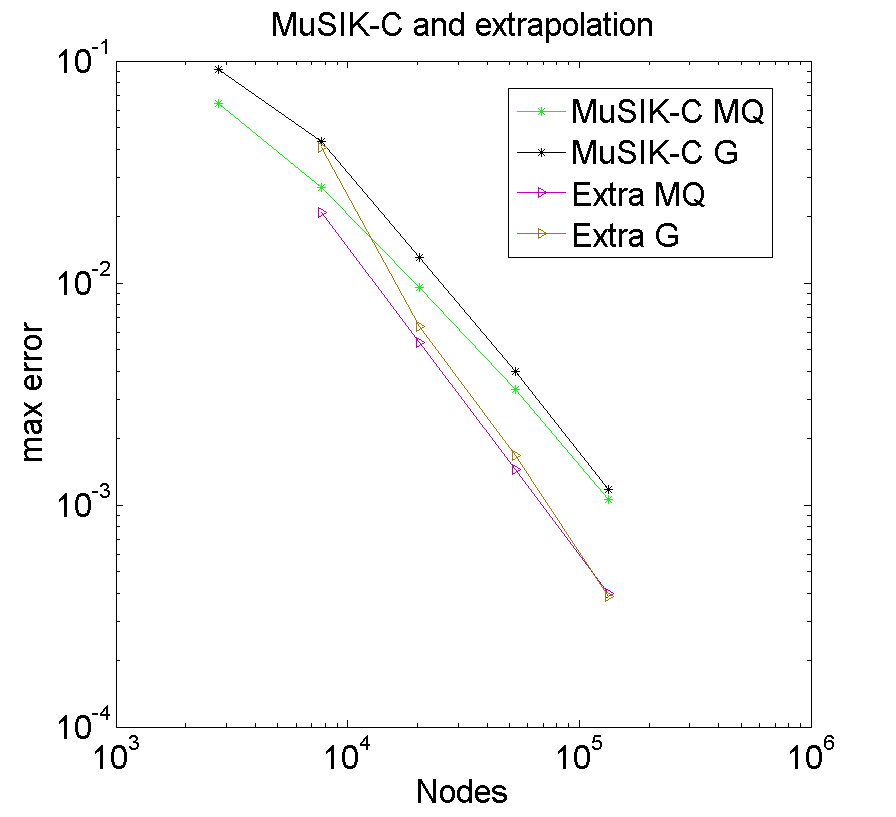}
\caption{MuSIK-C and correponding extrapolation with MQ and Gaussian when $C=2$ for Example \ref{PE 4D 2}.} \label{Fig:PE 4D 2}
\end{figure}

We see that the convergence rates for the non tensor product examples in Tables~\ref{Table PE 4D 2 C=2 MuSIK MQ} and \ref{Table PE 4D 2 C=2 MuSIK G} are worse than those for the tensor product case in Tables~\ref{Table PE 4D 1 C=2 MuSIK MQ} and \ref{Table PE 4D 1 C=2 MuSIK G}, for the range of numerical examples explored. Because we are not able to say if we are observing asymptotic rates, we cannot comment on whether or not the actual rates are worse.


\section{Conclusion}

The multilevel sparse grid kernel-based collocation (MuSIK-C) algorithm is used to solve elliptic and parabolic PDEs in up to four dimensions. In the  parabolic case we treat time as another space dimension and use the sparse grid on all dimensions. We use tensor product basis functions which are smooth but anisotropic, depending on the anisotropy in the sparse grid decomposision. We do numerical examples of both tensor product and non tensor product type. We compare our results to others in the literature.

Advantages of this method are that we reduce the overall complexity when compared to considering space and time separately. The use of smooth basis functions means that we have the possibility of spectral convergence orders, though the numerical results neither confirm nor deny this. MUSIK-C compares well in terms of convergence rate to the methods we compared with. In all methods improved convergence is observed with smoother basis functions.

In MuSIK-C the shape parameter of the smooth basis functions provides a smoothness parameter and we see an increase in condition number of the discrete systems as we increase the smoothness. It is the aim of future work to find pre-conditioning methods so that we can provide more numerically stable algorithms.

It is also the case that we can only use our method at this stage on domains which are simple to transform to hypercubes. This is a restriction when compared to other methods. However, this paper demonstrates that MuSIK-C has the potential to work in higher dimensions. The interpolation analogue MuSIK has been implemented in 10 dimensions and there is no reason why MuSIK-C cannot work in this dimension also.

Our numerical experiments indicate that MuSIK-C is more successful for tensor product problems, but that convergence is still observed for smooth non tensor product examples. We show that extrapolation can lead to improvements in error, though we do not achieve better convergence orders.

Future work will focus on solving ill-conditioning problems related to smoother basis functions, and on implementations in higher dimensions. 

\paragraph{Acknowledgements:} We are grateful to Peter Dong for running numerical experiments for comparison, and to Andrea Cangiani and Manolis Georgoulis for useful conversations related to this project. We also thank the referees for their careful reading and helpful comments; these have improved the presentation greatly (we hope).\\

%
%


\begin{thebibliography}{25}
\bibitem {Babuska}I. Babuska, F. Nobile and R. Tempone, \textit{A stochastic collocation method for elliptic partial differential equations with random input data}, SIAM J. Numer. Anal., 45 (2007), pp. 1005-1034.
\bibitem {Brown}D. Brown, L. Ling, E. Kansa, and J. Levesley, \textit{On approximate cardinal preconditioning methods for solving PDEs with radial basis functions}. Engineering Analysis with Boundary Elements, Vol. 29 (2005), pp. 343-353.
\bibitem {Bungartz}H.-J. Bungartz and M. Griebel, \textit{Sparse grids}, Acta Numer., 13 (2004), pp. 147-269.
\bibitem{Burden} R. L. Burden, J. D. Faires, and A. M. Burden, Numerical Analysis, 10th edition, Cengage Learning, 2016.

\bibitem {Delvos}F.-J. Delvos, \textit{d-variate Boolean interpolation}, J. Approx. Theory, 34 (1982), pp. 99-114.
\bibitem{Dong}
Z. Dong.
\newblock Personal communication.
\bibitem {Farrell}P. Farrell and H. Wendland, \textit{RBF multiscale collocation for second order elliptic boundary value problems}, SIAM J. Numer. Anal., Vol.51, No.4(2013), pp. 2403-2425.
\bibitem {Fasshauer1}G. E. Fasshauer, \textit{Solving partial differential equations by collocation with radial basis functions}, in Surface Fitting and Multiresolution Methods, A. Le Mehaute, C. Rabut, and L. L. Schumaker (eds.), Vanderbilt University Press, Nashville TN, 1997, 131-138.
\bibitem {Fasshauer}G. E. Fasshauer, \textit{Solving differential equations with radial basis functions: Multilevel methods and smoothing}, Adv. Comput. Math., 11 (1999), pp. 139-159.
\bibitem {Fornberg1}B. Fornberg and E. Lehto, \textit{Stabilization of RBF-generated finite difference methods for convective PDEs}, J. Comp. Phys., 230 (2011), pp. 2270-2285.

\bibitem {Fornberg}B. Fornberg and N. Flyer, \textit{Solving PDEs with radial basis function}, Acta Numerica, 24(2015), pp. 215-258.

\bibitem {Franke1}C. Franke and R. Schaback, \textit{Convergence order estimates of meshless collocation methods using radial basis functions}, Adv. Comput. Math., 8 (1998), pp. 381-399.

\bibitem {Franke2}C. Franke and R. Schaback, \textit{Solving partial differential equations by collocation using radial basis functions}, Appl. Math. \& Comp. 93(1998), pp. 73-82.

\bibitem {Ganapathysubramanian}B. Ganapathysubramanian and N. Zabaras, \textit{Sparse grid collocation schemes for stochastic natural convection problems}, Journal of Computational Physics, Vol. 225 (2007), pp. 652-685.

\bibitem {Garcke}J. Garcke and M. Griebel, \textit{On the parallelization of the sparse grid approach for data mining}, in Large-Scale Scientific Computations, Third International Conference, LSSC 2001.

\bibitem {Georgoulis}E. Georgoulis, J. Levesley, and F. Subhan, \textit{Multilevel sparse kernel-based interpolation}, SIAM J. Sci. Comput., 35 (2013), pp. 815-831.

\bibitem {Giesl}P. Giesl and H. Wendland, \textit{Meshless collocation: Error estimates with application to dynamical systems}, SIAM J. Numer. Anal., 45 (2007), pp. 1723-1741.

\bibitem {Griebel}M. Griebel, M. Schneider, and C. Zenger, \textit{A combination technique for the solution of sparse grid problems}, in Iterative methods in linear algebra (Brussels, 1991), North-Holland, Amsterdam, 1992, pp. 263-281.

\bibitem {Hales}S. J. Hales and J. Levesley, \textit{Error estimates for multilevel approximation using polyharmonic splines}, Numer. Algorithms, 30 (2002), pp. 1-10.

\bibitem {Hardy}R. L. Hardy, \textit{Multiquadric equations of topography and other irregular surfaces}, J. Geophys, Research, 45 (1971), pp. 1905-1915.

\bibitem {Hu}H. S. Hu and Z.-C. Li, \textit{Radial basis collocation methods for elliptic boundary value problems}, Comput. Math. Appl., 50 (2005), pp. 289-320.

\bibitem {Kansa1}E. J. Kansa, \textit{Multiquadrics: A scattered data approximation scheme with applications to computational
fluid-dynamics I. Surface approximations and partial derivative estimates}, Comput. Math. Appl., 19 (1990), pp. 127-145.

\bibitem {Kansa2}E. J. Kansa, \textit{Multiquadrics: A scattered data approximation scheme with applications to computational fluid-dynamics, part II: Solutions to parabolic, hyperbolic and elliptic partial differential equations}, Comput. Math. Appl., 19 (1990), pp. 147-161.

\bibitem {Kansa1986}E. J. Kansa, \textit{Application of Hardy’s multiquadric interpolation to hydrodynamics}, in Proc. 1986 Simul. Conf., Vol. 4, 1986, pp. 111-117.



\bibitem{Langer1}U. Langer, S. Moore, and M. Neumuller, \textit{Space-time isogeometric analysis of parabolic evolution equations}, Computer Methods in Applied Mechanics and Engineering, 306 (2016), pp. 342-263.

\bibitem {Micchelli}C. A. Micchelli, \textit{Interpolation of scattered data: distance matrices and conditionally positive definite functions}, Constructive Approximation , 2(1986), pp. 11-22.

\bibitem {Myers}D.E. Myers, S. De Iaco, D. Posa, L. De Cesare, \textit{Space-time radial basis functions}, Comput. Math. Appl., 43 (2002), pp. 539-549.

\bibitem {Nobile}F. Nobile, R. Tempone, and C. G. Webster, \textit{A sparse grid stochastic collocation method for partial differential equations with random input data}, SIAM J. Numer. Anal., 46 (2008), pp. 2309-2345.

\bibitem {Powell}M. J. D. Powell, \textit{The theory of radial basis function approximation in 1990, Advances in Numerical Analysis}, Vol.II, ed.W.Light, Oxford University Press, (1992), pp. 105-210.

\bibitem{rivlin} T. J. Rivlin, The Chebyshev Polynomials, John Wiley and Sons Inc., 1974.

\bibitem{SchabackU}  R. Schaback, \textit{Error estimates and condition numbers for radial basis function interpolation},  Adv. Comput. Math., 3 (1995), pp. 251–264.

\bibitem {Schaback}R. Schaback, \textit{Convergence of unsymmetric kernel‐based meshless collocation methods}, SIAM J. Numer. Anal., 45(1) (2007), pp. 333-351.

\bibitem{schwab} C. Schwab and R Stevenson, \textit{Space-time adaptive wavelet methods for parabolic
evolution problems}, Math. Comp., 78 (2009), pp. 1293–1318.

\bibitem {Shcherbakov}V. Shcherbakov, E. Larsson, \textit{Radial basis function partition of unity methods for pricing vanilla basket options}, Comput. Math. Appl., 71 (2016), pp. 185-200.


\bibitem {Usta}F. Usta, \textit{Sparse grid approximation with Gaussians} PhD, University of Leicester, 2015.

\bibitem {Wang1}Z. Wang, Q. Tang, W. Guo, and Y. Cheng, \textit{Sparse grid discontinuous Galerkin methods for high-dimensional elliptic equations}, Journal of Computational Physics, 314 (2016), pp. 244-263.

\end{thebibliography}
\end{document}